\documentclass[10pt,a4paper]{article}
\usepackage[utf8]{inputenc}	

\usepackage{amsfonts}
\usepackage{amssymb}			
\usepackage{amsmath}			
\usepackage{graphicx}			
\usepackage{hyperref}			
\usepackage[usenames,dvipsnames,svgnames,table]{xcolor}							
\hypersetup{colorlinks=true, citecolor=Blue, linkcolor=Blue, urlcolor=Blue}					
\usepackage{lineno}															

\usepackage[
maxnames=100,
url=false,
sorting=none,
firstinits=true,
backend=bibtex]{biblatex}
\newenvironment{bcases}{\left[} {\right\.}
\catcode`@=11
\def\caseswithdelim#1#2{\left#1\,\vcenter{\normalbaselines\m@th
  \ialign{\strut$##\hfil$&\quad##\hfil\crcr#2\crcr}}\right.}
\catcode`@=12
\usepackage[centering,total={17cm,26cm}]{geometry}								

\newtheorem{thm}{Theorem}

\newtheorem{defi}[thm]{Definition}

\newcommand{\R}{{\mathbb R}}

\DeclareMathOperator{\rem}{rem}
\DeclareMathOperator{\var}{var}
\DeclareMathOperator{\wind}{wind}
\DeclareMathOperator{\tr}{tr}

\DeclareMathOperator{\dsc}{dsc}
\DeclareMathOperator{\res}{res}

\begin{document}


\title{An algebro-geometric classification of spectral types of equilibria}

\author{Andrea Giacobbe\footnote{Università degli Studi di Catania, Dipartimento di Matematica e Informatica, Viale Andrea Doria 6, 95125 Catania, Italy, giacobbe@dmi.unict.it}}

\maketitle

\begin{abstract}
Given a vector field $X$ and an equilibrium $e$, we write three algebraic equations in the space of invariants of the linearisation of $X$ at $e$ which provide a geometric classification of all spectral types of the equilibrium. The loci defined by these equations correspond to definite types of bifurcations. The complement of such loci give a geometric decomposition of the space of invariants in open domains in which the equilibrium has a given spectral type. The usefulness of this approach lays in the fact that, when dealing with a parameter-dependent dynamical system, the pull-back of the loci from the space of invariants to the parameter space gives the bifurcation-decomposition of parameter space for the dynamical system at hand. In this article we apply the theory thoroughly in dimension 3 and 4, and we also give effective methods to explicitly compute the spectral indices.
\end{abstract}
\section{Introduction}

Assume to be given a vector field $X$ in a $m$ dimensional manifold $M$, and an equilibrium $e$ of $X$. Once chosen a local atlas $(x_1,..,x_m): U \subset M \to \R^m$, the linearization of $X$ at $e$ is the linear system
\begin{equation}\label{linear system}
\dot x = JX_{e} \, x, \qquad \text{where} \quad JX_{e} = \begin{pmatrix}
 \frac{\partial  X_1}{\partial x_1} (e) & \cdots  & \frac{\partial  X_1}{\partial x_m} (e) \\
 \vdots  & \ddots & \vdots  \\
 \frac{\partial  X_m}{\partial x_1} (e) & \cdots  & \frac{\partial  X_m}{\partial x_m} (e) \\
\end{pmatrix}.
\end{equation}
Among equilibria, the \emph{hyperbolic} equilibria (those for which the spectrum of $JX_{e}$ has no eigenvalues with zero real part) are particularly important, are generic, and the linear system \eqref{linear system} is topologically conjugate to the original one in a neighbourhood of the equilibrium \cite{1960.PAMS.Hartman, 1961.D.Grobman}. In these cases the Jordan blocks of the matrix $JX_{e}$ have a very mild effect on the quantitative form of solutions (secular terms), and no effect on the qualitative structure of solutions. It follows that hyperbolic equilibria can be classified using only the spectral decomposition of the matrix $JX_{e}$. In particular every equilibrium can be given an inertia-type decomposition using the names \emph{stable} and \emph{unstable} to indicate the sign of the eigenvalues and \emph{node} and \emph{focus} to indicate wether the eigenvalues are complex or real. 

In this article we classify hyperbolic equilibria using the symbols
\[
f_\beta^\alpha n_\delta^\gamma, \qquad \text{where} \quad \alpha,\beta,\gamma,\delta \in \mathbb N.
\]
With $f^\alpha_\beta$ we indicate the direct sum of $\alpha$ unstable foci and $\beta$ stable foci, with the symbol $n^\gamma_\delta$ we indicate the direct sum of $\gamma$ unstable nodes and $\delta$ stable nodes. Of course $2 \alpha + 2\beta + \gamma +\delta = m$, the dimension of the phase space. For the sake of clarity, in classical treaties the name stable node is typically referred to what we call stable double node $n_2$, the name unstable node to what we call unstable double node $n^2$, and the name saddle to what we call $n^1_1$, the direct product of a 1-dimensional stable and a 1-dimensional unstable node. We give the following definition.
\begin{defi}
Given a vector field $X$ and a hyperbolic equilibrium $e$, let
\begin{itemize}
\item $\alpha$ be the number of couples of complex conjugate eigenvalues of $JX_{e}$ with positive real part; 
\item $\beta$ be the number of couples of complex conjugate eigenvalues of $JX_{e}$ with negative real part; 
\item $\gamma$ be the number of positive real eigenvalues of $JX_{e}$;
\item $\delta$ be the number of negative real eigenvalues of $JX_{e}$.
\end{itemize}
We call the numbers $\alpha, \beta,\gamma,\delta$ \emph{spectral indices}, and we call the symbol $f_\beta^\alpha n_\delta^\gamma$ \emph{spectral type} of $e$.
\end{defi}

\begin{figure}
\begin{center}
\includegraphics[width=7cm]{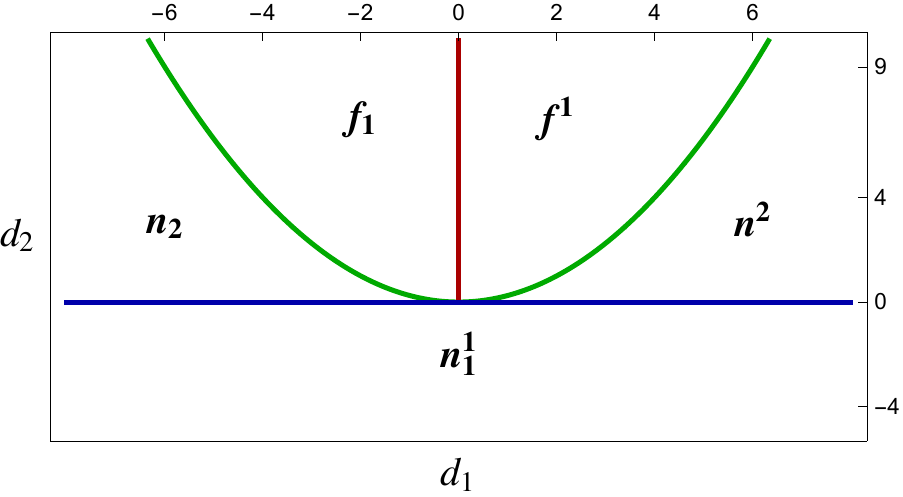}
\end{center}
\caption{The decomposition of the space of invariants $d_1$ (the trace) and $d_2$ (the determinant) associated to the equilibrium of a 2 degrees of freedom dynamical system. In each open set the spectral type of the equilibrium is specified.}\label{Marginal2}
\end{figure}

The investigation of the spectral type of an equilibrium in dimension 2 is trivial. In fact the linearisation of $X$ at $e$ yields a $2\times 2$ matrix whose characteristic polynomial is $p(\lambda) = \lambda^2 - d_1 \lambda + d_2$ where $d_1,d_2$ are the principal invariants of $JX_e$, that is $d_1 = \tr JX_{e}$ and $d_2 = \det JX_{e}$. The spectral type can be classified in the space of invariants in the well known diagram in Figure~\ref{Marginal2}.
In this diagram the red half-line is the semi-algebraic set $\mathcal{R} = \{(d_1,d_2) \,|\, d_1 = 0, d_2 > 0\}$, and it corresponds to the case in which $JX_{e}$ has a conjugate couple of non-zero purely imaginary eigenvalues. If we slightly move the principal invariants of $JX_{e}$ from right to left across such line, the sign of the real part of the two complex conjugate eigenvalues changes from positive to negative, meaning that the spectral type of the equilibrium changes from $f^1$ (unstable focus) to $f_1$ (stable focus). This event is called \emph{Hopf bifurcation}.

The blue line is the algebraic set $\mathcal{Z} = \{(d_1,d_2) \,|\,  d_2 = 0\}$, and it corresponds to the event of zero being an eigenvalue of $JX_{e}$. If the principal invariants of $JX_{e}$ slightly move across such line the sign of one real eigenvalue of $JX_{e}$ changes from positive to negative or vice-versa, meaning that the spectral type of the equilibrium changes from $n^2$ (unstable node) or $n_2$ (stable node) to $n_1^1$ (saddle). This event is called \emph{saddle-node bifurcation}.

The green curve is the algebraic set $\mathcal{D} = \{(d_1,d_2) \,|\, d_1^2 -4 d_2 = 0\}$, where $d_1^2 -4 d_2$ is the discriminant of the polynomial $p$, and it corresponds to the event of $p$ having a double root in the real line. If the principal invariants of $JX_{e}$ slightly move across such curve, two real eigenvalues with same sign become a double real eigenvalue with same sign and then transform into a couple of two complex conjugate eigenvalues with real part having the same sign of the double root (and hence same sign of the two distinct eigenvalues). This means that the spectral type of the equilibrium changes from $f^1$ to $n^2$ or from $f_1$ to $n_2$. This event is called \emph{focus-node bifurcation}.

The same type of analysis does not seem to have been carried out for higher dimensional systems, and in the literature it is possible to find only partial results, focussing on stability, which go under the name of Routh-Hurwitz conditions \cite{1895.MA.Hurwitz, 1877.Routh} (see also \cite{1971.JIMA.Barnett, 1991.CSSP.Anagnost.Desoer, 1992.CSM.Clark}).

In this article we set the analytical framework for a spectral analysis (Section~\ref{formal conditions}) and we perform a throughout investigation of the 3 and 4 dimensional cases, where the expressions are simple enough to be written and the results can be pictorially represented (Sections~\ref{3-dimensional case} and \ref{4-dimensional case}). The general expressions in higher dimensional cases become cumbersome (in Section~\ref{5 and 6 dimensional case} we show their aspect in dimension 5 and 6), and the geometric representation impossible. In Section~\ref{non-genericity} we discuss typical non-generic situation and their effect on our approach. In Section~\ref{7} we use Sturm's sequence and residue theorem to analytically compute the spectral indices. The general formulas, or even better the procedure to compute them, can be used in particular systems depending on few significant parameters, and gives a representation of the bifurcations of the system. To illustrate this fact we devote the last section.

\section{The formal conditions}\label{formal conditions}

Consider an equilibrium $e$ of a vector field $X$, and consider the characteristic polynomial of $JX_{e}$, which is the polynomial
\begin{equation}\label{car pol}
p(\lambda) = (-1)^m \lambda^m + (-1)^{m-1} d_1 \lambda^{m-1} + \cdots - d_{m-1} \lambda +  d_m.
\end{equation}
where $d_1 = \tr JX_e,...,d_m = \det JX_e$ are the principal invariants of the matrix $JX_e$. To every choice of principal invariants $(d_1,....,d_m)$ there corresponds a polynomial $p$ (the characteristic polynomial) and a family of roots of $p$ (the eigenvalues), that is a definite spectrum of $JX_{e}$, that is a definite spectral type of the equilibrium $e$. We will abuse the terminology and refer to \emph{the spectral type of the point $(d_1,...,d_m)$} in the space of invariants.

A change in spectral type of $e$ can take place only when very specific events take place. These events define algebraic varieties that are stratified varieties. Such algebraic varieties decompose the space of invariants in domains. We call \emph{marginal} all the points $(d_1,...,d_m)$ of the space of invariants at which a change of spectral type is taking place or, in other words, any point at the boundary of domains whose points have a given spectral type. The generic elementary changes in the spectral type for the equilibrium $e$ are all and only one of the following:
\begin{itemize}
\item[(z)] a single real root changes sign, which corresponds to the change of spectral type
\[
f^\alpha_\beta n^{\gamma+1}_\delta  \leftrightarrow f^\alpha_\beta n^\gamma_{\delta+1};
\]
\item[(r)] the real part of two complex conjugate roots changes sign, which corresponds to the change of spectral type
\[
f^{\alpha+1}_\beta n^\gamma_\delta  \leftrightarrow f^\alpha_{\beta+1} n^\gamma_\delta;
\]

\item[(d)] two complex conjugate roots collide in the real axis on a non-zero real number, and then separate into two real roots or viceversa, which corresponds to one of the two possible changes of spectral type, depending on the sign of the real part of the roots 
\[
f^{\alpha+1}_\beta n^\gamma_\delta \leftrightarrow  f^\alpha_\beta n^{\gamma+2}_\delta \qquad \text{or} \qquad f^\alpha_{\beta+1} n^\gamma_\delta \leftrightarrow f^\alpha_\beta n^\gamma_{\delta+2}.
\]
\end{itemize}

Condition (d) should be divided in two conditions: (d$^+$), corresponding to bifurcation $f^{\alpha+1}_\beta n^\gamma_\delta \leftrightarrow  f^\alpha_\beta n^{\gamma+2}_\delta$ and (d$^-$), corresponding to bifurcation $f^\alpha_{\beta+1} n^\gamma_\delta \leftrightarrow f^\alpha_\beta n^\gamma_{\delta+2}$. We will briefly address this issue at the end of this section. We prefer to keep this analysis out of the picture for clarity.

In each of these situation a very specific event must take place. The case (z) is particularly simple to treat. At the marginality corresponding to a bifurcation in which a root changes sign, zero must be a root of the characteristic polynomial $p$, and hence the function $\zeta(d_1,...,d_m) = d_m$ must vanish. We denote
\[
\mathcal Z = \{(d_1,...,d_m) \in \R^m \,|\, \zeta(d_1,...,d_m)= 0\}.
\]
This condition gives a hyperplane in invariant space which separates domains in which the spectral type of the equilibrium changes according to (z).

Case (d) is more involved. This type of bifurcation takes place when the characteristic polynomial has a double real root, which can happen only if the function $\delta(d_1,...,d_m) = \dsc(p)$, the discriminant of the polynomial $p$, vanishes. In this case the relevant algebraic variety $\mathcal D$ is a subvariety of
\[
\widetilde{\mathcal D} = \{(d_1,...,d_m) \in \mathbb R^m \,|\, \delta(d_1,...,d_m)  = 0\}.
\]
The reason for being a subvariety is due to the fact that multiple roots of $p$ could be outside the real axis. Such spurious solutions are strata of the variety $\widetilde{\mathcal D}$ which have higher codimension, and can be easily distinguished from the true marginal points (they can be so easily distinguished that they are often overseen). We will see in the case $n= 4$ that the marginal variety $\mathcal D$ differs form $\widetilde{\mathcal D}$ for a 1 dimensional curve which is the analogous of the thread emanating from the swallowtail singularity. One can formally define the marginal region $\mathcal D$ as the closure of $\widetilde{\mathcal D}^{m-1}$, where $\widetilde{\mathcal D}^{m-1}$ is the union of all $m-1$ dimensional strata of $\widetilde{\mathcal D}$ (see \cite{2012.Arnold.Gusein-Zade.Varchenko.1} for a discussion on stratifications).

The bifurcation of type (r) is the most complicate to treat. If we are at such marginality then the characteristic polynomial $p$ has two conjugate, purely imaginary roots, that we indicate $i \mu$, $-i \mu$, with $\mu$ real and non-zero. Let us denote
\begin{equation}\label{pr pi}
\begin{cases}
p^r(\mu) = d_m - d_{m-2} \mu^2 + d_{m-4} \mu^4 + \cdots = \sum_{j = 0}^{\lfloor m/2 \rfloor} (-1)^j d_{m-2j } \mu^{2j} \\[5pt]
p^i(\mu) = d_{m-1} - d_{m-3} \mu^2 + d_{m-5} \mu^4 + \cdots = \sum_{j = 0}^{\lfloor m/2 \rfloor} (-1)^j d_{m-1-2j} \mu^{2j}.
\end{cases}
\end{equation}
the two polynomials such that $p(i \mu) = p^r(\mu) - i \mu p^i(\mu)$ (in these expressions we agree that $d_0 = 1$, and $\lfloor m/2 \rfloor$ indicates the integer part of $m/2$). These two polynomials have degrees
\[
\deg (p^r) = 
\left[\begin{matrix}
m &\text{ if } m \text{ is even}\\
m-1  &\text{ if } m \text{ is odd},
\end{matrix}\right.
\qquad 
\deg(p^i) = \left[\begin{matrix}
m-2  & \text{ if } m \text{ is even}\\
m-1  & \text{ if } m \text{ is odd}.
\end{matrix}\right.
\]

At marginal points of type (r) the two polynomials $p^r$ and $p^i$  must have two common real roots $\pm \mu$. Unfortunately, both polynomials are in the variable $\mu^2$, and hence a codimension one condition is that these two polynomials have two common real solution or that they have two complex conjugate purely imaginary solutions. The polynomials $p^r$, $p^i$ have common roots precisely when the function $\widetilde\rho(d_1,...,d_m) = \res(p^r,p^i)$, the resultant of the two polynomials, vanishes. The above mentioned fact implies that, generically, it is not the entire variety
\[
\widetilde{\mathcal R} = \{ (d_1,...,d_m) \in \mathbb R^m \,|\, \widetilde\rho(d_1,....,d_m) = 0\}
\]
which corresponds to the marginality at exam, but only the semialgebraic variety that corresponds to a common double \emph{real} root of the polynomials $p^r$, $p^i$. 

We can rephrase the considerations above using the two polynomials $q^r(\nu)$ and $q^i(\nu)$ such that $p^r(\mu) = q^r(\mu^2)$ and $p^i(\mu) = q^i(\mu^2)$. The two polynomials are
\begin{equation}\label{qr qi}
q^r(\nu) = \sum_{j = 0}^{\lfloor m/2 \rfloor} (-1)^j d_{m-2j} \nu^j, \qquad q^i(\nu) =  \sum_{j = 0}^{\lfloor m/2 \rfloor} (-1)^j d_{m-1-2j} \nu^j
\end{equation}
and their degrees are
\[
\deg (q^r) = 
\left[\begin{matrix}
\frac m2 &\text{ if } m \text{ is even}\\[5pt]
\frac{m-1}2  &\text{ if } m \text{ is odd}
\end{matrix}\right.
\qquad 
\deg(q^i) = \left[\begin{matrix}
\frac m2 - 1  & \text{ if } m \text{ is even}\\[5pt]
\frac{m-1}2  & \text{ if } m \text{ is odd.}
\end{matrix}\right.
\]
The bifurcation of type (r) takes place when these two polynomials have a common \emph{positive} real root. 

Once again, the two polynomials have common roots when the function  $\rho(d_1,...,d_m) = \res(q^r,q^i)$ vanishes (observe that $\widetilde\rho = \rho^2$). But the vanishing of $\rho$ corresponds to the condition that the two polynomials have a common real root, while the marginal locus we are looking for is the subvariety of $\widetilde{\mathcal R}$ that corresponds to points $(d_1,...,d_m)$ whose associated polynomials $q^r$, $q^i$ have a \emph{common positive real root}. We hence need an extra condition to ensure that the common root is positive.

This condition can be obtained using Euclid's division algorithm. In fact the ultimate remainder of Euclid's division algorithm applied to $q^r$ and $q^i$ is a degree zero polynomial (a real number) that is the resultant $\rho$. When the resultant is zero, the penultimate remainder of the Euclid's division algorithm applied to $q^r$ and $q^i$ is a degree 1 polynomial whose only root $\sigma(d_1,...,d_m)$ is the common root that must exist given that the resultant is zero. This fact, which holds under generic assumptions, is what solves our dilemma, and we can state that
\[
\mathcal R = \{ (d_1,...,d_m) \in \mathbb R^m \,|\, \rho(d_1,...,d_m) = 0, \sigma(d_1,...,d_m) > 0\}.
\]

The remainders of Euclid's division of two polynomials have a fundamental role in the investigation of roots of polynomials, and they generate the so called \emph{Sylvester sequence} \cite{2013.Gondim.deMoralesMelo.Russo}. We summarise the discussion above in a definition and a main theorem.

\begin{defi}
We call \emph{determinant locus} the set $\mathcal Z$, \emph{discriminant locus} the set $\mathcal D$ and \emph{resultant locus} the set $\mathcal R$. We call \emph{marginal locus} the union of the three loci. We call \emph{marginal points} the points of the marginal locus.
\end{defi}

\begin{thm}\label{main}
Given a vector field $X$ and an equilibrium $e$ of $X$. The marginal locus in invariant space is the union of three algebraic varieties: $\mathcal Z$, $\mathcal R$, and $\mathcal D$. These varieties decompose the space of invariants in domains with a specific spectral type. Across the marginal locus the spectral type of $e$ varies according to a precise rule depending on the three possibilities (z), (r), (d) listed above.

The variety $\mathcal Z$ is the hyperplane $\{\zeta = 0\}$ with $\zeta$ the determinant of $JX_e$. With possible lower-dimensional artefacts, the algebraic variety $\mathcal D$ is the variety $\{\delta = 0\}$ with $\delta$ the discriminant of $p$, the characteristic polynomial of $JX_e$; the semialgebraic variety $\mathcal R$ is the semialgebraic variety $\{\rho = 0, \sigma > 0\}$ with $\rho$ the resultant and $\sigma$ the unique root of the penultimate Euclid's remainder of the two polynomials $q^r,q^i$ defined in \eqref{qr qi}.
\end{thm}

We conclude observing that the argument of the penultimate Euclid's remainder applied to $p$ and $p'$ does give information on what type of bifurcation is taking place between $f^{\alpha+1}_\beta n^\gamma_\delta \leftrightarrow  f^\alpha_\beta n^{\gamma+2}_\delta$ and $f^\alpha_{\beta+1} n^\gamma_\delta \leftrightarrow f^\alpha_\beta n^\gamma_{\delta+2}$. In fact the last remainder of Euclid's division algorithm applied to $p$, $p'$ is the discriminant. When the discriminant is zero, the root $\tau(d_1,...,d_m)$ of the penultimate remainder (a degree one polynomial) is the real double root of $p$, and hence its sign will discriminate between the two possible bifurcations: the unstable focus$\leftrightarrow$unstable node (d$^+$) or the stable focus$\leftrightarrow$stable node (d$^-$).

\section{The 3-dimensional case}\label{3-dimensional case}

Let us use Theorem~\ref{main} to classify the spectral type of hyperbolic equilibria of a 3-dimensional system. Consider $p(\lambda) = - \lambda^3 + d_1 \lambda^2 - d_2 \lambda + d_3$, the characteristic polynomial of a $3\times 3$ matrix, where $d_i$ are the invariants of the matrix, i.e.\ $d_3$ is the determinant, $d_1$ is the trace, $d_2$ is the sum of the determinants of the three principal $2\times 2$ minors.

\begin{figure}[h]
\begin{center}
\includegraphics[width=5cm]{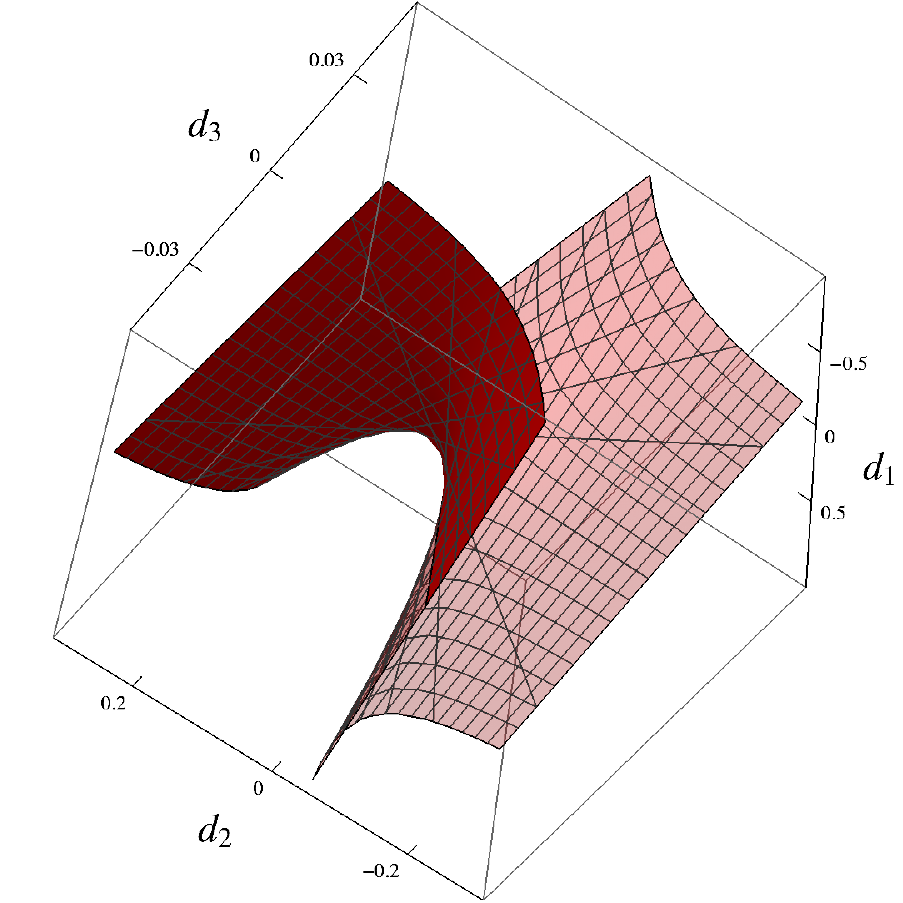}
\includegraphics[width=5cm]{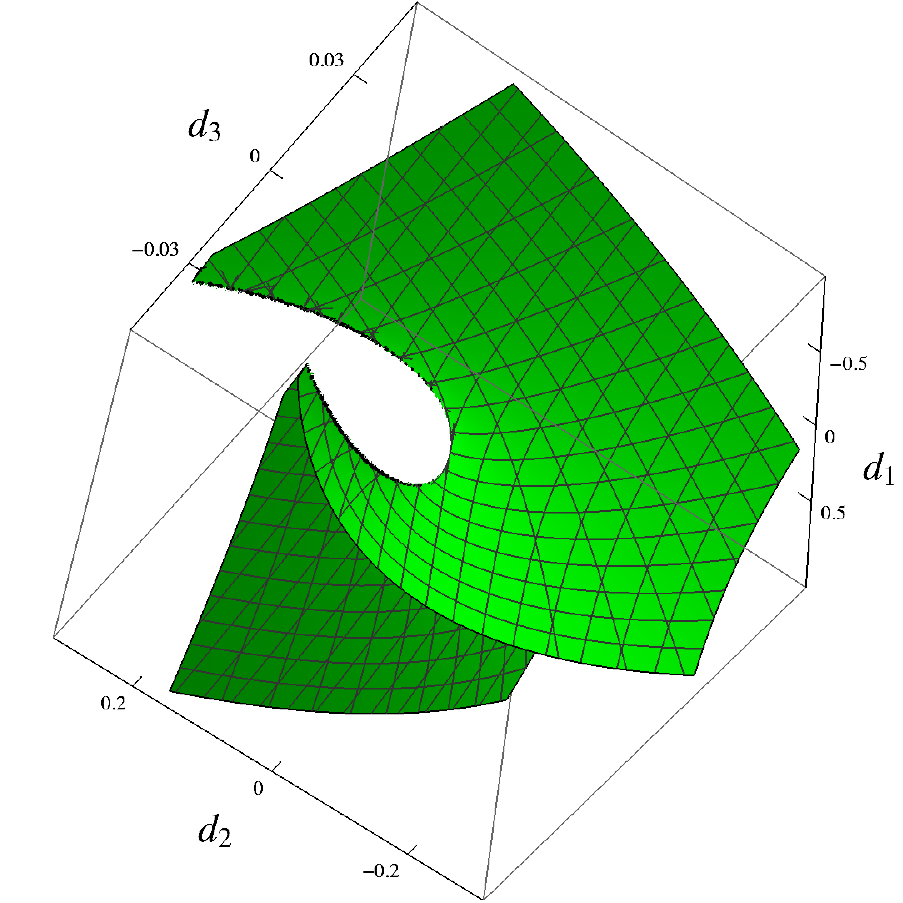}
\includegraphics[width=5cm]{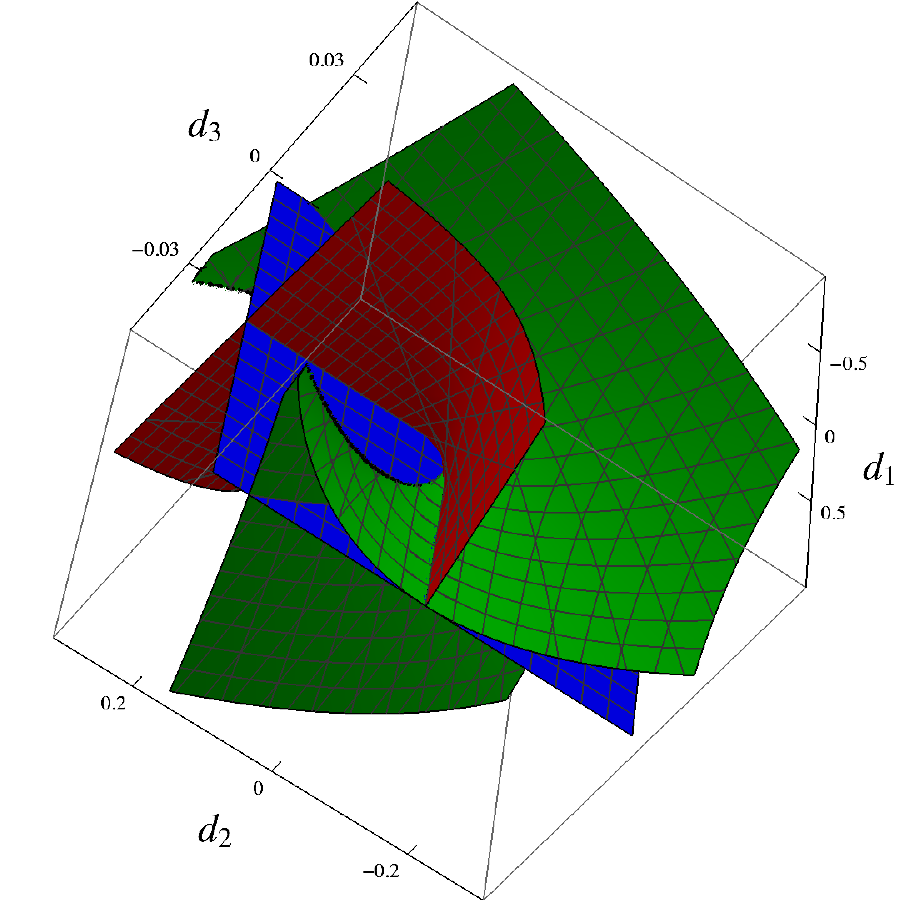}
\end{center}
\caption{In the left pane the resultant locus $\mathcal R$ (solid red),  a semialgebraic subvariety of the hyperbolic paraboloid $\widetilde{\mathcal R}$ (union of solid and transparent red). In the center pane the discriminant locus $\mathcal D$. In the right pane the two varieties together with the determinant locus $\mathcal Z$. The complement of the marginal varieties correspond to domains in which the spectral type does not change.}
\label{n=3}
\end{figure}

The hyperplane $\mathcal Z = \{(d_1,d_2,d_3) \,|\, d_3 = 0\}$ is easily drawn. Also the discriminant of the characteristic polynomial is easy to compute, and is
\[
\delta = -4 d_3 d_1^3+d_2^2 d_1^2+18 d_2 d_3 d_1-4 d_2^3-27 d_3^2.
\]
The corresponding discriminant locus $\mathcal D = \{(d_1,d_2,d_3) \in \mathbb R^3  \,|\, \delta(d_1,d_2,d_3) = 0\}$ is drawn in Figure~\ref{n=3} center pane. In this low dimensional case the polynomials $p$ cannot have a double complex root, since this event can take place only when the polynomial has degree at least four, hence in this case $\widetilde{\mathcal D} = \mathcal D$. The interesting feature of this algebraic set is that $\mathcal D$ displays a line of cusp points corresponding to a triple root in the real axis.

For the variety $\mathcal R$ we must consider the two polynomials $q^r = -d_1 \nu + d_3$ and $q^i = - \nu + d_2$. They have common positive real roots only if $d_3/d_1 = d_2$ and $d_2 > 0$. The resultant of the two polynomials is in fact $\rho = d_3 -d_1 d_2$. In this case the penultimate remainder is $q^i$ itself (the system is very low-dimensional). It follows that $\sigma = d_2$. In Figure~\ref{n=3} left pane a picture of the resultant locus $\mathcal R$, (in transparent red is represented $\widetilde{\mathcal R} \setminus \mathcal R$) in the right pane a cumulative picture of the three loci $\mathcal Z \cup \mathcal D \cup \mathcal R$.

\begin{figure}[h]
\begin{center}
\includegraphics[width=5cm]{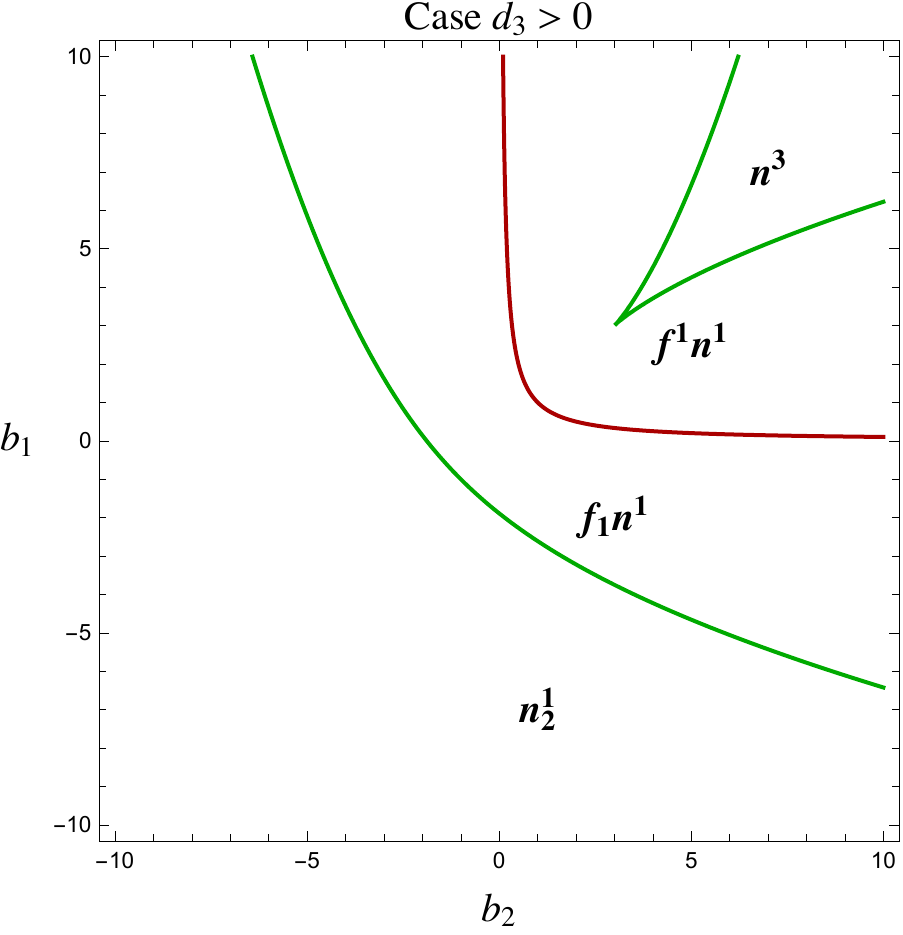}
\includegraphics[width=5cm]{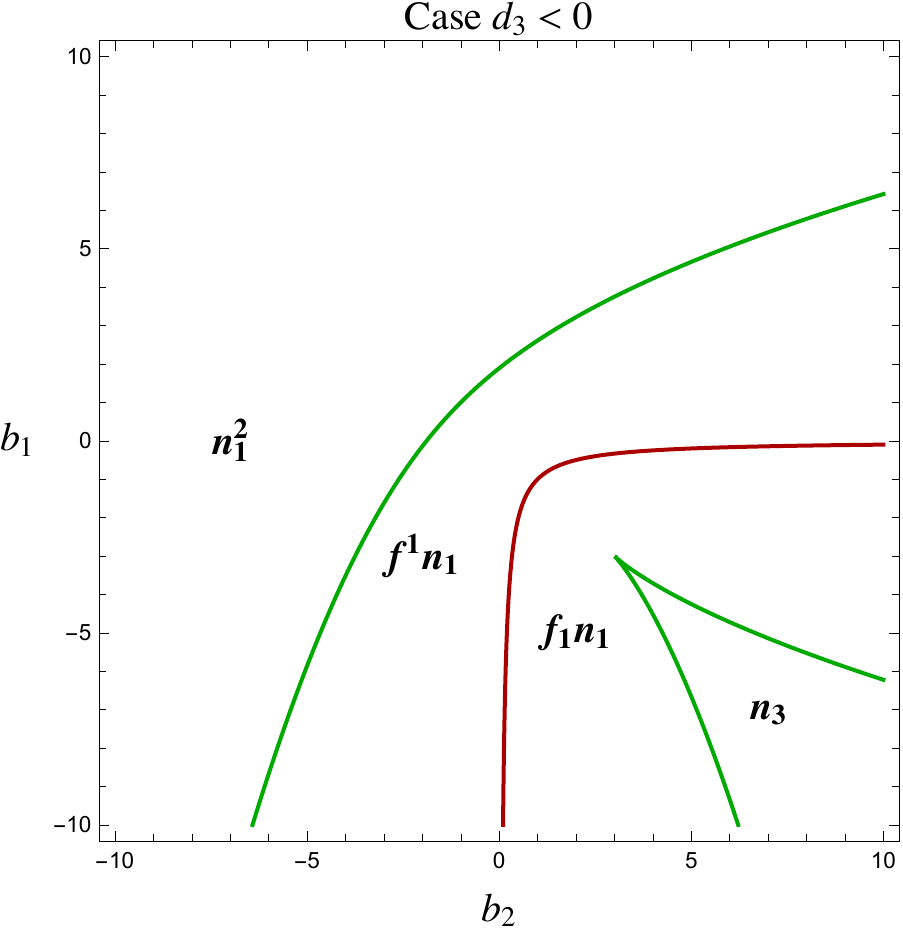}
\end{center}
\caption{Tomographics sections of Figure~\ref{n=3} with specification of the spectral type of the equilibrium in the complement of the marginal loci.}
\label{n=3 tomography}
\end{figure}

Although complete, the above pictures lacks the indication of spectral type in each of the connected component complement of the marginal loci. For this reason we use a simple consideration to choose a representative slice of Figure~\ref{n=3}. A positive rescaling of the vector field $X$ amounts to a rescaling of the spectral parameter $\lambda$ (if the rescaling is non positive, there will be an exchange of the stable indices $\beta$, $\delta$ with the unstable ones $\alpha$, $\gamma$). Assume therefore that $d_3 \not = 0$. We can positively rescale the variable $\lambda$ by posing $\lambda = k \widetilde\lambda$ with $k \in \mathbb R^+$ and obtain the polynomial
\[
\widetilde p(\widetilde \lambda) = k^3 \left(\widetilde \lambda^3 - \frac{d_1}k  \widetilde \lambda^2 + \frac{d_2}{k^2} \widetilde \lambda - \frac{d_3}{k^3}\right)
\]
whose roots have same spectral type of the roots of $p$. We hence shall choose $k = \sqrt[3]{|d_3|}$ and investigate two possibilities:
\[
\left[
\begin{matrix}
p^- = \widetilde \lambda^3 - b_1 \widetilde \lambda^2 + b_2 \widetilde \lambda - 1 & \text{ if } d_3 < 0\\[3pt]
p^+ = \widetilde \lambda^3 - b_1 \widetilde \lambda^2 + b_2 \widetilde \lambda + 1 & \text{ if } d_3 > 0
\end{matrix}
\right.
\]
where $b_1 = d_1 /\sqrt[3]{|d_3|}$ and $b_2 = d_2 /\sqrt[3]{|d_3|}$. The discriminants of these two polynomials are the already computed discriminant $\delta$ with the substitutions $d_3 \to \pm 1$, $d_2 \to b_2$, $d_1 \to b_1$, that is
\[
\delta^\pm = b_1^2 b_2^2  - 4 b_1^3 \mp 4 b_2^3  \pm 18 b_1 b_2  - 27.
\]
The vanishing of $\delta^\pm$ always corresponds to a double real root of $p$ except at the codimension 2 cusp point, which corresponds to a triple real root (see the green curves of Figure~\ref{n=3 tomography}). Such variety is the marginal state separating the case of $p$ possessing three real roots to the case of $p$ possessing two complex conjugate roots and one real root.

After substitution the resultant is $\rho^\pm = \pm 1 - b_1b_2$, and of course its relevant submanifold is the one in which the polynomials $q^r = - b_1 \nu \pm 1 $ and $q^i = - \nu + b_2$ have a common positive real root. It follows that $b_2 > 0$. In Figure~\ref{n=3 tomography} we show how the marginal loci separate the space of invariants in regions of homogeneous spectral type, and how across each locus the change in spectral type is determined by the locus being crossed.

The sign of $\tau(d_1,d_2,d_3) = (d_1 d_2-9 d_3)/(2 \left(d_1^2-3 d_2\right))$ separates $\mathcal D$ in two marginal loci, $\mathcal D^+$ across which two positive real roots become a couple of complex conjugate roots with positive real part and $\mathcal D^-$ across which two negative real roots become a couple of complex conjugate roots with negative real part.

\section{The 4-dimensional case}\label{4-dimensional case}

In the 4-dimensional case the characteristic polynomial of a $4\times 4$ matrix is $p(\lambda) = \lambda^4 - d_1 \lambda^3 + d_2 \lambda^2 - d_3 \lambda + d_4$. Also in this case the coefficients $d_1,d_2,d_3,d_4$ are the principal invariants of the matrix. Also in this case zero is an eigenvalue if and only if $d_4 = 0$. So that the determinant locus $\mathcal Z$ is a hyperplane whose points separate regions in which one eigenvalue changes sign. The discriminant of $p$ is
\begin{multline*}
\delta = -27 d_4^2 d_1^4-4 d_3^3 d_1^3+18 d_2 d_3 d_4 d_1^3+d_2^2 d_3^2 d_1^2+144 d_2 d_4^2 d_1^2 - 4  d_2^3 d_4 d_1^2-6 d_3^2 d_4 d_1^2 + \\
+ 18 d_2 d_3^3 d_1 -192 d_3 d_4^2 d_1 - 80 d_2^2 d_3 d_4  d_1-27 d_3^4+256 d_4^3-4 d_2^3 d_3^2-128 d_2^2 d_4^2+16 d_2^4 d_4+144 d_2 d_3^2 d_4,
\end{multline*}
while the two polynomials needed to define $\mathcal R$ are $q^r =\nu ^2 -d_2 \nu + d_4 $, $q^i = d_1 \nu - d_3$, from which it follows that
\[
\rho = d_4 d_1^2-d_2 d_3 d_1+d_3^2, \qquad \sigma = d_1 d_3.
\] 

These functions are what is needed to investigate bifurcations of 4-dimensional systems, but in this case the space of invariants is 4-dimensional. 

We proceed as done in the 3-dimensional case and, assuming $d_4 \not = 0$ we apply a positive rescaling of the vector field, which amounts to a positive rescaling of the variable $\lambda$, hence reducing the parameters down to three. Pose $\lambda = k \widetilde \lambda$ with $k  \in \mathbb R^+$, the polynomial $p$ becomes
\[
\widetilde p(\widetilde \lambda) = k^4 \left(\widetilde \lambda^4 - \frac{d_1}k \widetilde \lambda^3 + \frac{d_2}{k^2} \widetilde \lambda^2 - \frac{d_3}{k^3} \widetilde \lambda+ \frac{d_4}{k^4}\right),
\]
whose roots have same spectral type of the roots of $p$. We hence shall choose $k = \sqrt[4]{|d_4|}$ and investigate two possibilities:
\[
\left[
\begin{matrix}
p^- = \widetilde \lambda^4 - b_1 \widetilde \lambda^3 + b_2 \widetilde \lambda^2 - b_3 \widetilde \lambda - 1 & \text{ if } d_4 < 0\\[3pt]
p^+ = \widetilde \lambda^4 - b_1 \widetilde \lambda^3 + b_2 \widetilde \lambda^2 - b_3 \widetilde \lambda+ 1& \text{ if } d_4 > 0.
\end{matrix}
\right.
\]
(also in this case $b_j = d_j/\sqrt[4]{|d_4|}$ for $j = 1,2,3$.) The discriminant of these two polynomials is
\begin{multline}
\delta^\pm = -27 b_1^4+18 b_1^3 b_2 b_3 - 4 b_1^3  b_3^3 - 4 b_1^2 b_2^3 + b_1^2 b_2^2   b_3^2  \pm 144 b_1^2 b_2 \mp 6 b_1^2 b_3^2 \mp 80b_1b_2^2 b_3 \pm 18 b_1 b_2 b_3^3-192 b_1 b_3+ \\
\pm 16 b_2^4 \mp 4 b_2^3 b_3^2 - 128 b_2^2  + 144 b_2 b_3^2-27 b_3^4 \pm 256.
\end{multline}
In this case the variety $\delta^+ = 0$ does posses a thread, corresponding to the codimension 2 degeneracy associated to the coincidence of two couples of complex conjugate solutions. These degeneracies do not correspond to marginal points, but they are of fundamental interest, being origin of interesting monodromic effects. The variety $\mathcal D$ has also codimension 2 strata which are cusp singularities, corresponding to triple real roots, which are edges to regular strata of codimension 1 corresponding to double real roots, and it also has two point singularities that are swallotails, and correspond to a quadruple real root \cite{2012.Arnold.Gusein-Zade.Varchenko.1}.

With this reduction from the space of invariants $(d_1,d_2,d_3,d_4)$ to the space $(b_1,b_2,b_3)$, the resultant becomes $\rho^\pm = b_1^2 - b_1 b_2 b_3 \pm b_3^2$ and, when the resultant is zero, the common root of $q^r$ and $q^i$ is $\sigma^\pm = b_1/b_3$. This fact can be deduced from the abstract approach of Section~\ref{formal conditions}, but it can also be obtained by direct computation from the two polynomials $q^r = \nu^2 - b_2 \nu  \pm 1$ and $q^i = b_3 \nu - b_1$.

\begin{figure}[ht]
\begin{center}
\includegraphics[width=5cm]{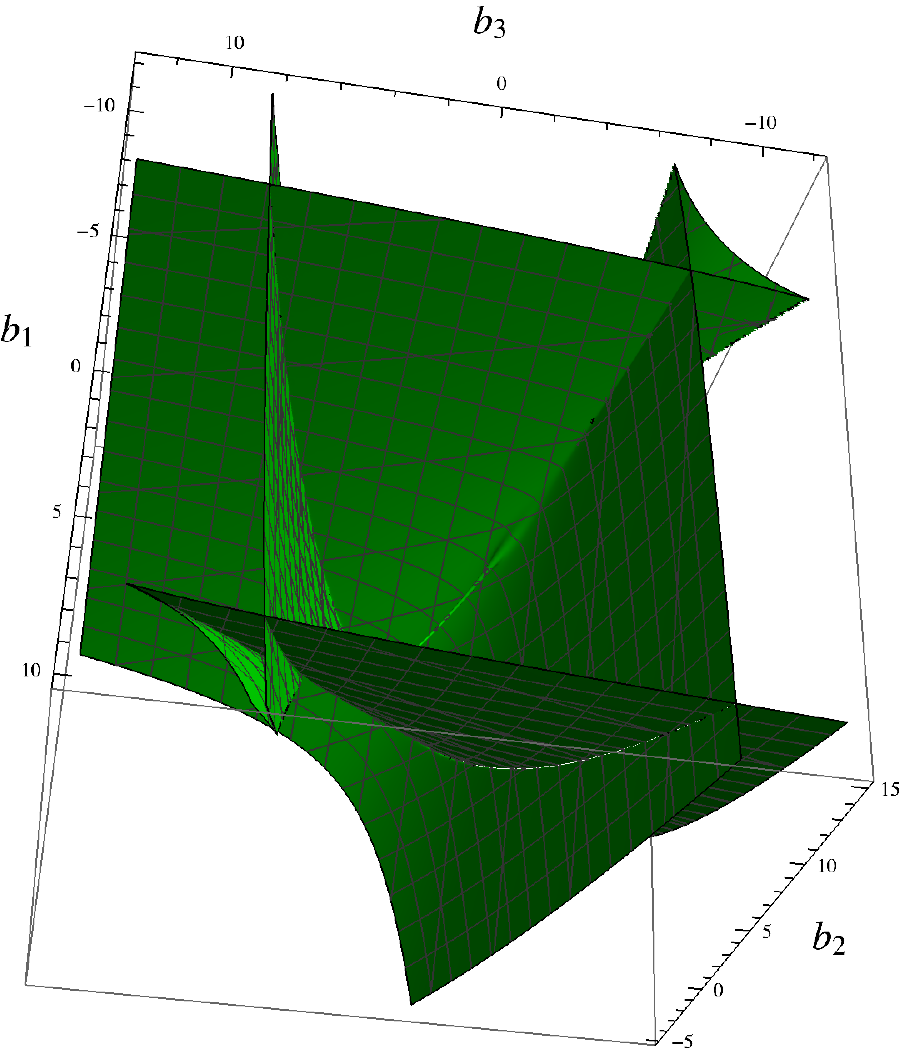}
\includegraphics[width=5cm]{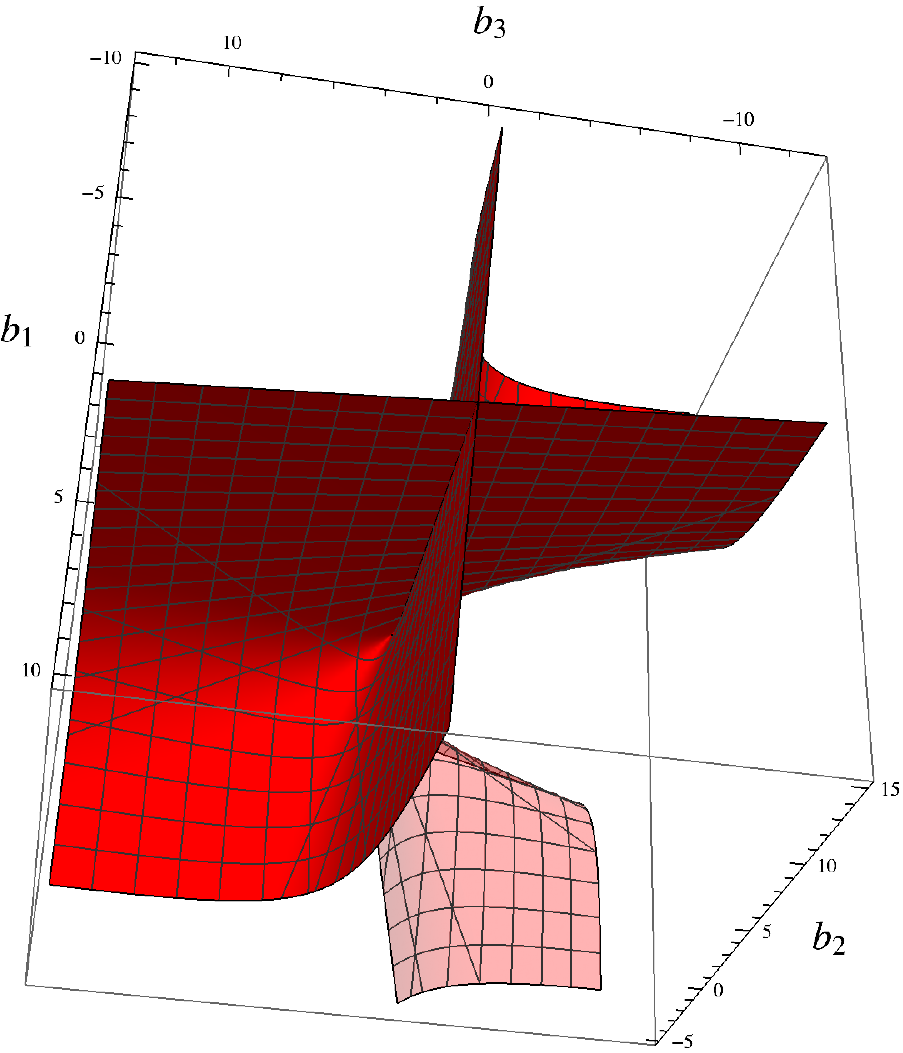}
\includegraphics[width=5cm]{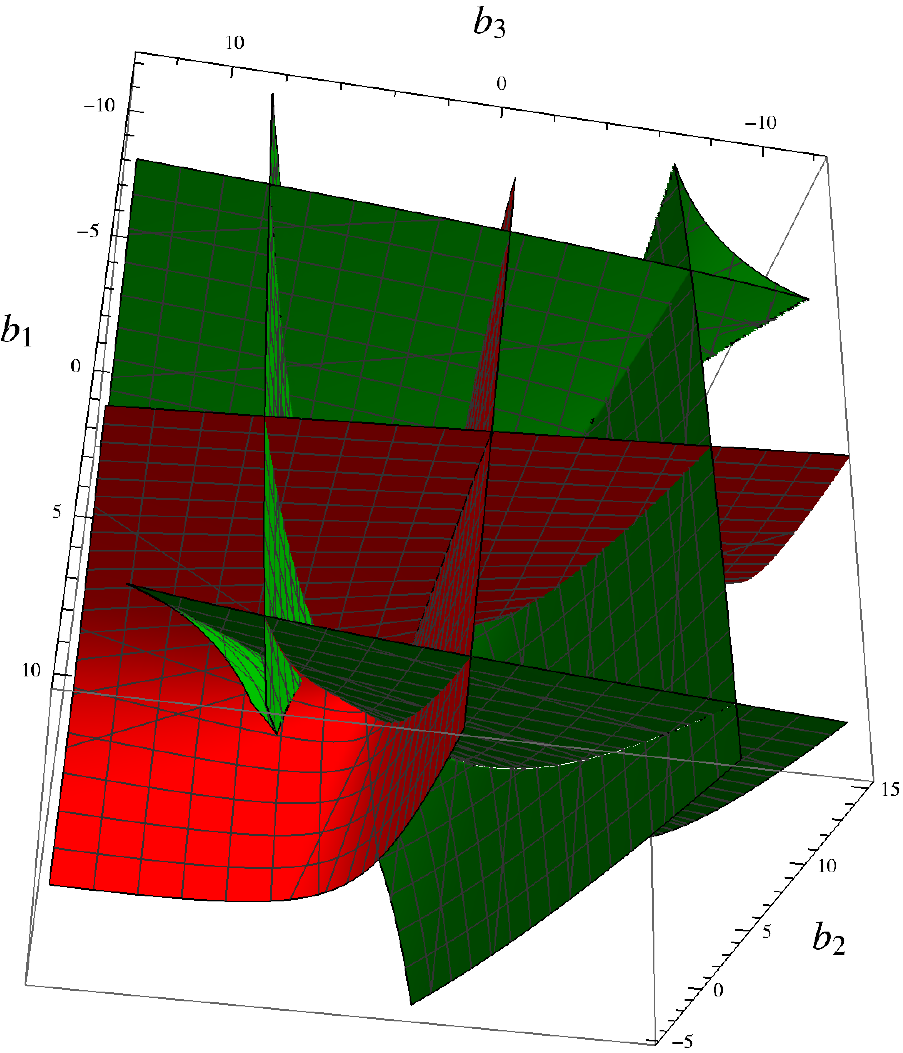}
\end{center}
\caption{Case in which $d_4 >0$. In the left pane the marginal locus $\widetilde{\mathcal D}$, which differs from $\mathcal D$ only for the 1-dimensional thread that connects the two 0-dimensional strata (sligtly visible at the center of the image). In the center pane the locus $\widetilde{\mathcal R}$ in red (solid and transparent) and $\mathcal R$ in solid red. In the right pane $\mathcal D$ and $\mathcal R$ together.}\label{n=4 d4>0}
\end{figure}

\begin{figure}[h]
\begin{center}
\includegraphics[width=5cm]{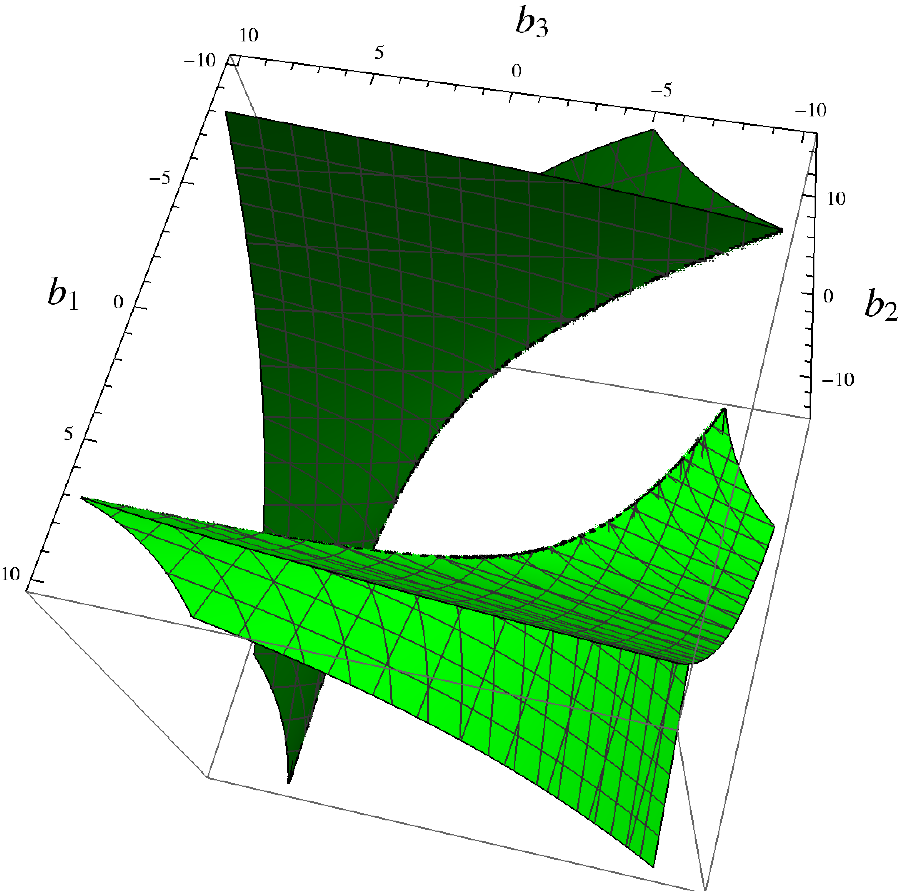}
\includegraphics[width=5cm]{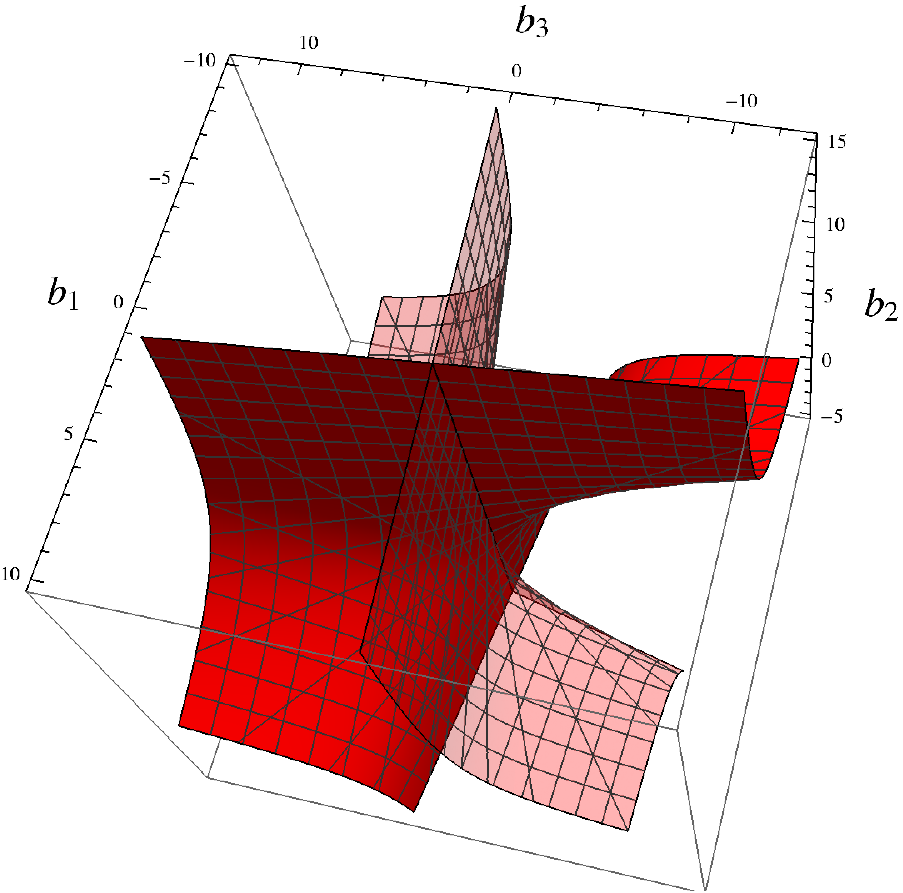}
\includegraphics[width=5cm]{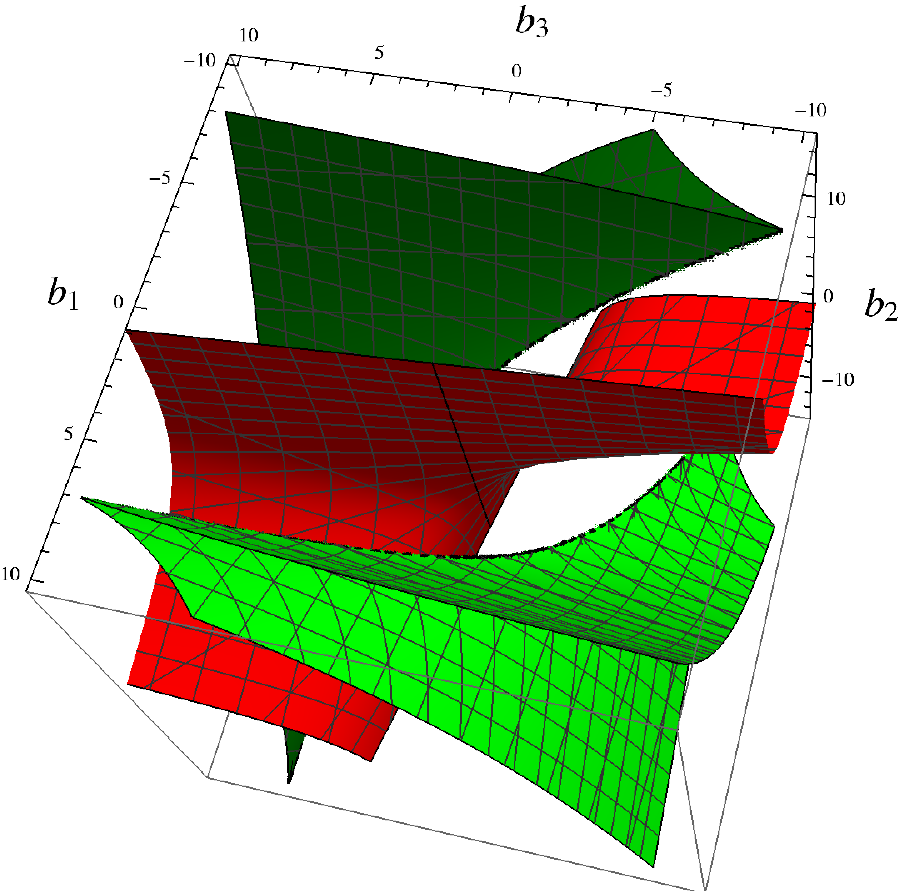}
\end{center}
\caption{Case in which $d_4 < 0$. In the left pane the marginal locus $\mathcal D$. In the center pane the locus $\widetilde{\mathcal R}$ in red (solid and transparent) and $\mathcal R$ in solid red. In the right pane $\mathcal D$ and $\mathcal R$ together.}\label{n=4 d4<0}
\end{figure}

Representations of the loci when $d_4 > 0$ are given in Figure~\ref{n=4 d4>0}, while representations of the loci when $d_4 < 0$ are given in Figure~\ref{n=4 d4<0}.

\begin{figure}[ht]
\begin{center}
\includegraphics[width=5cm]{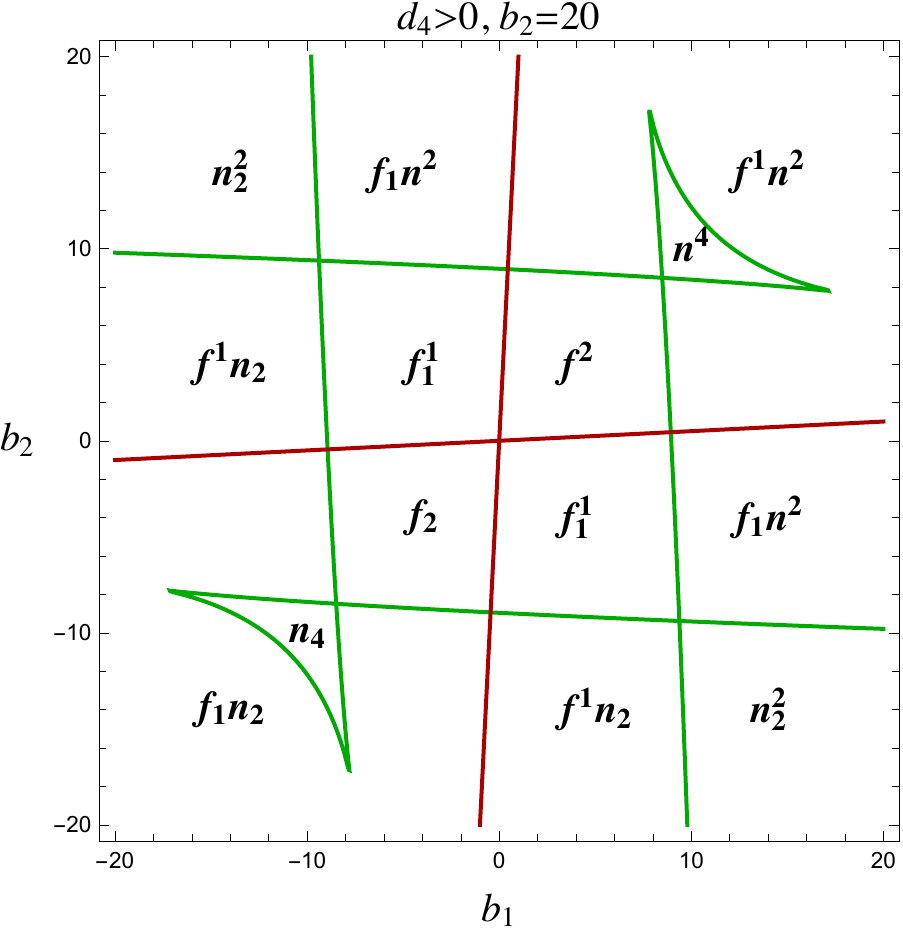}
\includegraphics[width=5cm]{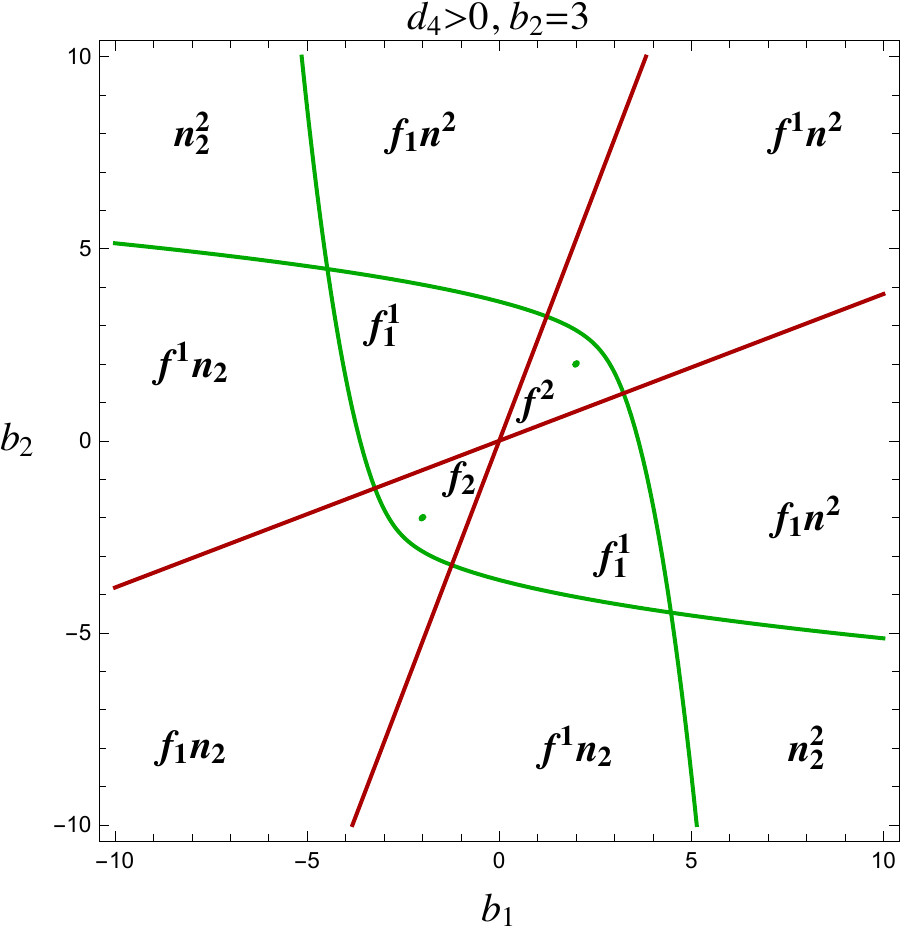}
\includegraphics[width=5cm]{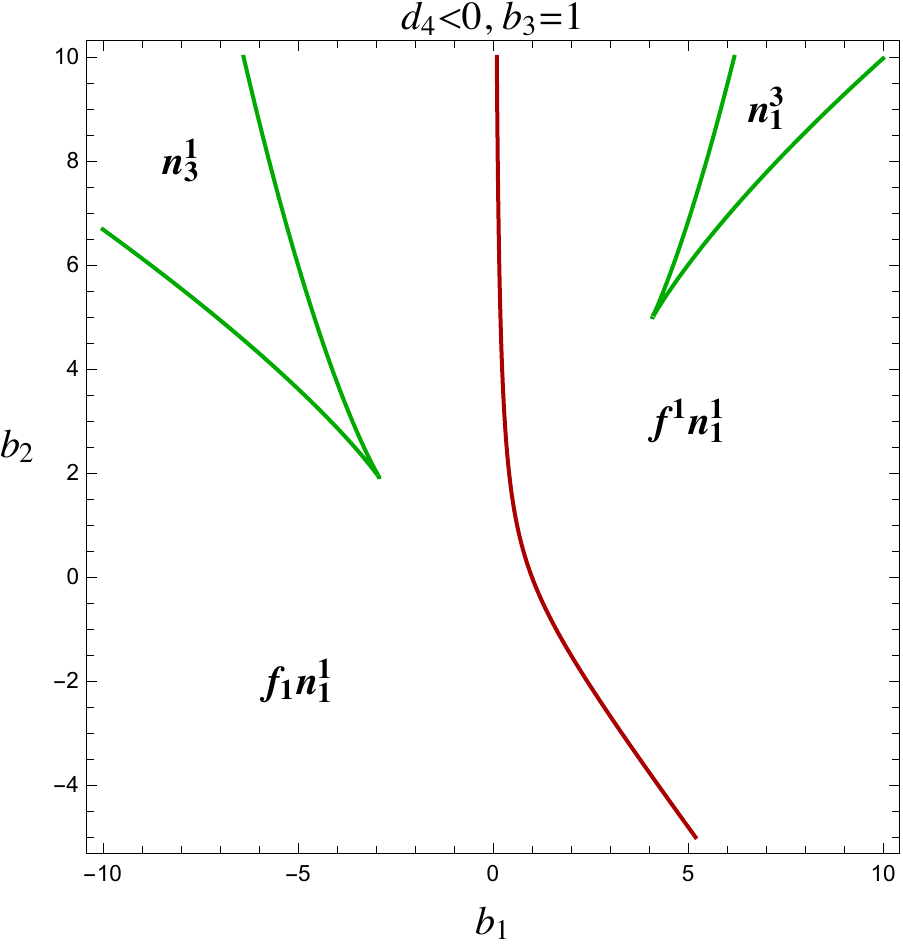}
\end{center}
\caption{In the left and center panes, tomographic sections of the 3-dimensional Figure~\ref{n=4 d4>0} at two different $b_2$ levels. The labels of the different domains give an exhaustive description of the spectral types that can be met. In the center pane are visible two isolated points (drawn as infinitesimal circles for technical reasons). These points are the intersection of the tomographic plane with the thread in which a focus-focus degeneracy takes place. In the right pane, a section of the 3-dimensional Figure~\ref{n=4 d4<0} at a chosen $b_3$ level. The labels in the different domains give an exhaustive description of the spectral types that can be met.}
\label{tomography n=4}
\end{figure}

We complete this investigation by making a tomography at a few representative level sets and indicating the spectral type of the connected components of the complement of the loci. In Figure~\ref{tomography n=4} we use the same convention made for the 3-dimensional case, that is red line for the resultant locus (across which a change from stable focus to unstable focus takes place) and green line for the discriminant locus (across which a change focus to node takes place, preserving the type of stability).

While in the odd-dimensional case ($m$ odd) the regions with $d_m >0$ are equivariantly symmetric with those having $d_m < 0$, which means that the loci are the same, and the spectral types change switching the upper indices with the lower ones, in the even-dimensional case such equivariance is lost, because the vector field $-X$ has linearisation with the same determinant of the linearisation of $X$. This means that in the even dimensional case the figures with $d_m$ positive or negative have a $\mathbb Z_2$ symmetry that the figures in the odd dimensional case do not possess. In invariant space such $\mathbb Z_2$ symmetry can be explicitly written and is 
\[
\left[
\begin{matrix}
(d_1,d_2,d_3,...,d_{n-1}, d_m) \to (-d_1,d_2,-d_3,...-d_{n-1},d_m)  &\text{ if } m \text{ is even}\\[5pt]
(d_1,d_2,d_3,...,d_{n-1}, d_m) \to (-d_1,d_2,-d_3,...d_{n-1},- d_m) &\text{ if } m \text{ is odd}.
\end{matrix}
\right.
\]
This symmetry becomes an equivariance for the odd-dimensional case between points with positive and negative determinant, clearly visible in Figure~\ref{n=3} and in Figure~\ref{n=3 tomography} between left and right pane, while in the even-dimensional case it becomes the invariance $(b_1,b_2,b_3,...,b_{n-1}) \to (-b_1,b_2,-b_3,...-b_{n-1})$ clearly visible in Figures~\ref{n=4 d4>0}, \ref{n=4 d4<0}, and in the first two panes of Figure~\ref{tomography n=4}, where it corresponds to the transformation $(b_1,b_2,b_3) \to (-b_1,b_2,-b_3)$.

Another $\mathbb Z_2$ symmetry for the reduced case, that is $(b_1,...,b_{n-1})$-space, is visible both in the odd and even dimensional case. It is related to the substitution of the spectral parameter $\lambda = 1/\widetilde\lambda$. This symmetry corresponds to an alternating sign inversion of the parameters 
\[
\left[
\begin{matrix}
(b_1,...,b_{n-1}) \to (\mp b_{n-1}, \pm b_{n-2}, ..., \pm b_2, \mp b_1) &\text{ in the even case}\\[5pt]
(b_1,...,b_{n-1}) \to (\mp b_{n-1}, \pm b_{n-2}, ..., \mp b_2, \pm b_1) &\text{ in the odd case.}
\end{matrix}
\right.
\]

\section{The 5 and 6-dimensional cases}\label{5 and 6 dimensional case}

In the 5-dimensional case the polynomial $p$ and the two polynomials $q^r,q^i$ are
\[
p = -\lambda ^5 + d_1 \lambda ^4 - d_2 \lambda ^3 + d_3 \lambda ^2 - d_4 \lambda + d_5, \qquad q^r = d_1 \nu ^2-d_3 \nu +d_5, \qquad q^i = -\nu ^2 + d_2 \nu - d_4.
\]
It follows that
\begin{multline*}
\delta = 256 d_5^3 d_1^5-27 d_4^4 d_1^4-128 d_3^2 d_5^2 d_1^4-192 d_2 d_4 d_5^2 d_1^4+144 d_3 d_4^2  d_5 d_1^4+18 d_2 d_3 d_4^3 d_1^3-1600 d_2 d_5^3 d_1^3 +\\
-4 d_3^3 d_4^2 d_1^3+144 d_2^2 d_3 d_5^2 d_1^3+160 d_3 d_4 d_5^2 d_1^3+16 d_3^4 d_5 d_1^3-36 d_4^3 d_5 d_1^3-6 d_2^2 d_4^2 d_5 d_1^3-80 d_2 d_3^2 d_4 d_5 d_1^3+144 d_2 d_4^4 d_1^2 + \\
-4 d_2^3 d_4^3 d_1^2-6 d_3^2 d_4^3 d_1^2+2000 d_3 d_5^3 d_1^2+d_2^2 d_3^2 d_4^2 d_1^2-27 d_2^4 d_5^2 d_1^2+560 d_2  d_3^2 d_5^2 d_1^2-50 d_4^2 d_5^2 d_1^2+1020 d_2^2 d_4 d_5^2 d_1^2-4 d_2^2 d_3^3 d_5 d_1^2+ \\
-746 d_2 d_3 d_4^2 d_5 d_1^2+24 d_3^3 d_4 d_5 d_1^2+18 d_2^3 d_3 d_4 d_5 d_1^2-192   d_3 d_4^4 d_1-80 d_2^2 d_3 d_4^3 d_1+2250 d_2^2 d_5^3 d_1-2500 d_4 d_5^3 d_1+18 d_2 d_3^3 d_4^2 d_1+ \\
 -900 d_3^3 d_5^2 d_1-630 d_2^3 d_3 d_5^2 d_1-2050 d_2 d_3 d_4 d_5^2  d_1-72 d_2 d_3^4 d_5 d_1+160 d_2 d_4^3 d_5 d_1+24 d_2^3 d_4^2 d_5 d_1+1020 d_3^2 d_4^2   d_5 d_1+356 d_2^2 d_3^2 d_4 d_5 d_1+\\
 + 256 d_4^5-128 d_2^2 d_4^4+3125 d_5^4+16 d_2^4 d_4^3+144 d_2 d_3^2 d_4^3-3750 d_2 d_3 d_5^3-27 d_3^4 d_4^2-4 d_2^3 d_3^2 d_4^2+108  d_2^5 d_5^2+825 d_2^2 d_3^2 d_5^2+\\
 + 2000 d_2 d_4^2 d_5^2-900 d_2^3 d_4 d_5^2+2250 d_3^2 d_4 d_5^2+108 d_3^5 d_5+16 d_2^3 d_3^3 d_5-1600 d_3 d_4^3 d_5+560 d_2^2 d_3 d_4^2 d_5-630 d_2 d_3^3 d_4 d_5-72 d_2^4 d_3 d_4 d_5,
 \end{multline*}
while
\[
\rho = d_1 d_5 d_2^2-d_1 d_3 d_4 d_2-d_3 d_5 d_2+d_1^2 d_4^2+d_5^2+d_3^2 d_4-2 d_1 d_4 d_5.
\]
The penultimate remainder of Euclid's division algorithm applied to $q^r$, $q^i$ is $\left(d_1 d_2-d_3\right) \nu -d_1 d_4+d_5$, and hence its root is $\sigma = (d_1 d_4-d_5)/(d_1 d_2-d_3)$. The positivity of this function is equivalent to the positivity of the polynomial function
\[
\sigma = d_2 d_4 d_1^2-d_3 d_4 d_1-d_2 d_5 d_1+d_3 d_5,
\]
which is preferable to use being a polynomial expression. In the 6-dimensional case the polynomial $p$ and the two polynomials $q^r,q^i$ are
\[
p = \lambda^6 - d_1 \lambda ^5 + d_2 \lambda ^4 - d_3 \lambda^3 + d_4 \lambda^2 - d_5 \lambda + d_6, \qquad q^r = -\nu ^3 + d_2 \nu ^2-d_4 \nu + d_6, \qquad q^i = -d_1 \nu ^2+d_3 \nu -d_5.
\]
We spare the reader from seeing the expression of $\delta$, while 
\begin{multline*}
\rho = -d_6^2 d_1^3-d_4^2 d_5 d_1^2+d_3 d_4 d_6 d_1^2+2 d_2 d_5 d_6 d_1^2-d_2^2 d_5^2 d_1+2 d_4 d_5^2 d_1+\\
+ d_2 d_3 d_4 d_5 d_1-d_2 d_3^2 d_6 d_1-3 d_3 d_5 d_6 d_1-d_5^3+d_2 d_3 d_5^2-d_3^2 d_4 d_5+d_3^3 d_6.
 \end{multline*}
In this case the penultimate remainder of Euclid's division algorithm applied to $q^r$, $q^i$ is the polynomial
\[
\left(-\frac{d_3^2}{d_1^2}+\frac{d_2 d_3}{d_1}-d_4+\frac{d_5}{d_1}\right) \nu -\frac{d_2  d_5}{d_1}+\frac{d_3 d_5}{d_1^2}+d_6.
\]
It follows that $\sigma$ can be chosen as
\[
\sigma = \left(d_4 d_1^2 - d_1 d_2 d_3 - d_1 d_5+d_3^2\right) \left(d_6 d_1^2-d_2 d_5 d_1+d_3 d_5\right).
\]

\section{Non-genericity}\label{non-genericity}
\begin{figure}[ht]
\begin{center}
\includegraphics[width=5cm]{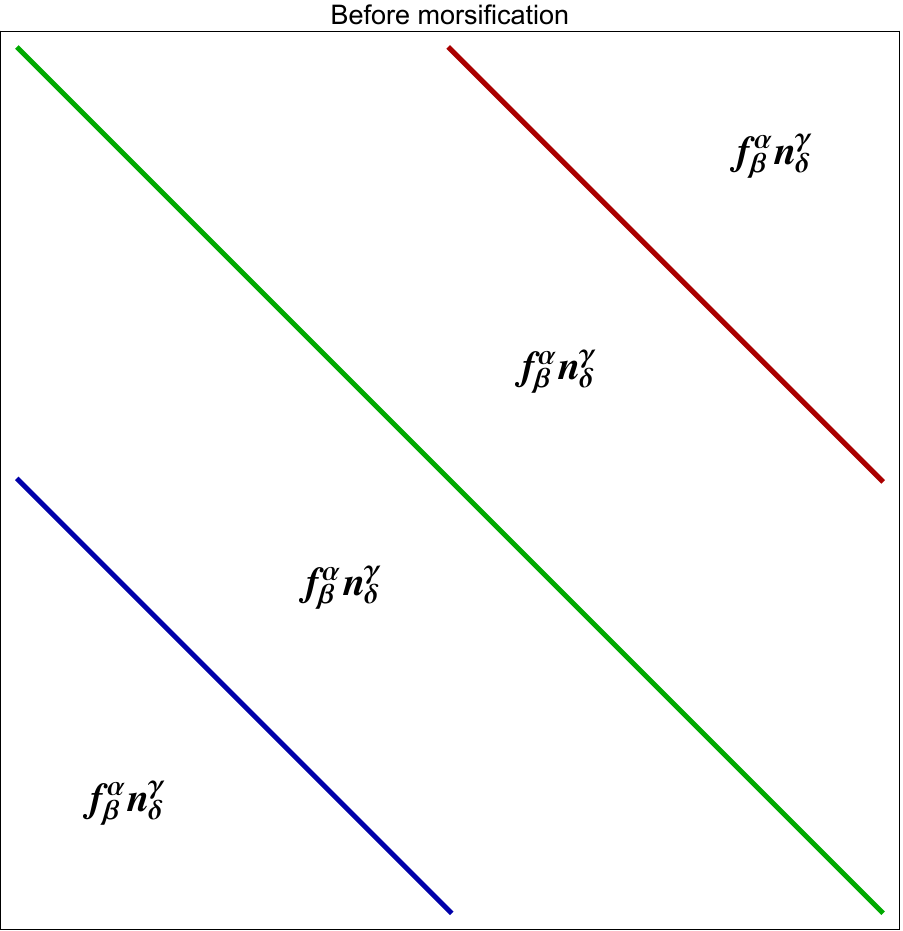} \qquad \includegraphics[width=5cm]{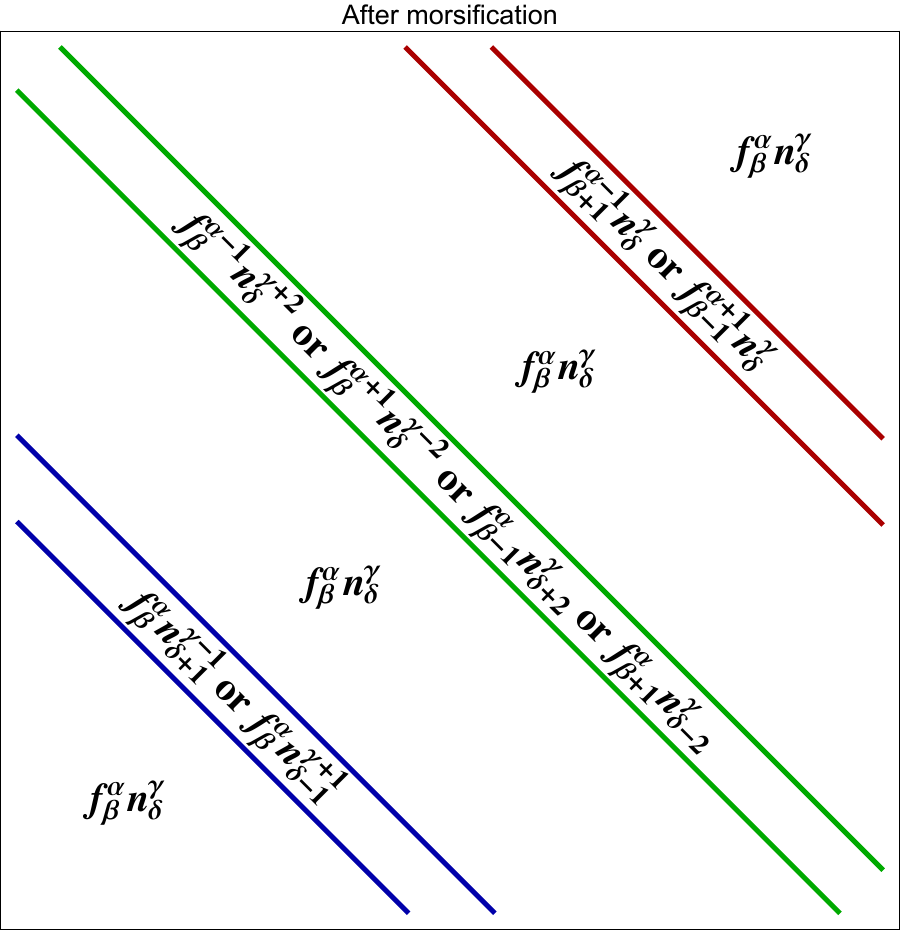}
\end{center}
\caption{These 2 panes represent the effect of morsification on the higher dimensional component of marginal loci ($\mathcal Z$, $\mathcal D$, $\mathcal R$ respectively) when, across them, the sign of their defining function ($\zeta$, $\delta$, $\rho$ respectively) does not change sign.}\label{morsification}
\end{figure}

Some words must be spent on the hypothesis we made that the bifurcations are generic. In many cases the equilibrium of a parameter-dependent dynamical system undergoes bifurcation at some parameters, but the bifurcation does not evolve with a change of spectral type. The three typical cases that can take place are:
\begin{itemize}
\item[(z$_{deg}$)] one real eigenvalue becomes zero and then moves back to the same real semi-axis from which it came from;
\item[(r$_{deg}$)]
two complex conjugate eigenvalues touch the purely imaginary axis and then move back to the same half-plane from which they came from;
\item[(d$_{deg}$)]
two real eigenvalues collide in the real axis but they then separate once again in the real axis instead of separating in a couple of complex-conjugate egenvalues (or the same event with two complex conjugate that separate back to two complex conjugate).
\end{itemize}
These three events are non-generic. Non generic events such as these (and other more subtle) take often place because the vector field $X$ has some symmetry, and the principal invariants of the matrix $JX_e$ are not free to move in an open domain of invariant space. From the point of view of the present treatment this fact can be fully understood by applying a small generic perturbation to the vector field $X$. In singularity theory this is called \emph{morsification}. An example of what this would mean is represented in Figure~\ref{morsification}.

It is possible that in a chosen vector field our approach will display some of these degenerate features (and even more dramatic ones). The non-generic events described above can be easily determined: in all such cases the algebraic function ($\zeta$, $\delta$, or $\rho$) whose zeroes are the marginal locus being crossed do not change sign across the marginal variety in parameter space. Observe that a software for numerical representation of level-sets will typically not draw such varieties, since a contour hypersurface is typically detected by tracing the changes in sign of the defining function. We will see this fact in the application of Section~\ref{application}.

One last word must be spent on a remarkably interesting family of dynamical systems: the Hamiltonian vector fields. Hamiltonian vector fields are non-generic, and the characteristic polynomial of their linearisation at the equilibria are always polynomials in $\lambda^2$. This implies that the degree of the characteristic polynomial is always even, and that $p^i$ is always zero. Lots can be said on such vector fields, but this is outside the purpose of this article. For such vector fields our approach can give results only after non-Hamiltonian morsification. We will dedicate to Hamiltonian vector fields an ad-hoc investigation.

\section{Determination of the indices}\label{7}

In the previous section we made clear a fact that was not pointed at explicitely in the first part of the manuscript: at each marginal locus one type of bifurcation must take place, but there is uncertainty on which direction this bifurcations is taking place (from positive to negative or vice-versa, from two real to two complex conjugate or vice-versa). This is never unclear at low dimensions, and the best way to settle the question in applications is by checking at  appropriate points of the connected components complement of the marginal loci $\mathcal Z, \mathcal D, \mathcal R$. Nonetheless, a theoretical discussion requires a non-numerical determination of the precise spectral type, meaning that we would like to express the spectral type only as a fuction of the invariants $d_1,...,d_m$. This can be done almost explicitly using Sturm's sequences and residue theorem.

Let us begin by considering the indices $\gamma$ and $\delta$. They indicate the number of positive real zeroes and of negative real zeroes of the characteristic polynomial $p$. This information can be easily extracted from Sturm's theorem \cite{1835.Sturm}. Consider $p_0 = p$ and $p_1 = p'$, and define for $j= 2,..., m$ the polynomial $p_j = - \rem(p_{j-2},p_{j-1})$ where $\rem(p_{j-2},p_{j-1})$ is the remainder of the Euclidean division of $p_{j-2}$ by $p_{j-1}$. Let $a_1,a_2,....,a_m$ be the list of constant terms of the polynomials $p_j$ and let $b_1,....,b_m$ be the coefficients of the top power of the polynomials $p_j$.  We have that
\[
\gamma = \var(a_1,a_2,....,a_m) - \var(b_1,....,b_m), \delta =  \var(a_1,a_2,....,a_m) - \var((-1)^{1-m}b_1,....,(-1)^{j-m}b_j,...,(-1)^{m-m}b_m)
\]
where with $\var$ we mean the number of changes in sign of consecutive elements of the list (we should say "once the zeroes have been removed", but such event is not generic). This could appear a mysterious formula, but it can be expressed explicitly in terms of the invariants.

The integer numbers $2\alpha+\gamma$ and $2\beta+\delta$ can be computed using the residue theorem. In fact if the polynomial $p$ has simple zeroes, the function $p'/p$ has simple poles with residue 1 at every zero of $p$ \cite{1974.Henrici}. It follows that the number of zeros that $p$ possesses in a region $\Omega$ is $(2\pi i)^{-1} \int_{\partial \Omega} p'(z)/p(z) dz$ (with the usual convention on orientation of boundaries). By the lemma of the big circle, one has that
\[
2 \alpha + \gamma = \frac 1{2 \pi i} \lim_{R \to + \infty} \left(\int_R^{-R} \frac{p'(i s)}{p(i s)} i ds + \int_{\Gamma_R} \frac{p'(z)}{p(z)} dz\right) \quad 2 \beta + \delta = \frac 1{2 \pi i} \lim_{R \to + \infty} \left(\int_{-R}^R \frac{p'(i s)}{p(i s)} i ds + \int_{\Delta_R} \frac{p'(z)}{p(z)} dz\right),
\]
where $\Gamma_R$ is the big half-circle parametrised by $\Gamma_R(\vartheta) = R(\cos\vartheta + i \sin\vartheta)$, $\vartheta \in [-\pi/2,\pi/2]$ and $\Delta_R$ is the other big half-circle parametrised by $\Delta_R(\vartheta) = R(\cos\vartheta + i \sin\vartheta)$, $\vartheta \in [\pi/2,3\pi/2]$ . Concentrating on $2\alpha+\gamma$ we have that the integral along the big half-circle $\Gamma_R$ gives $n \pi i $, from which it follows that
\begin{multline*}
2 \alpha + \gamma - \frac n2  = \frac 1{2 \pi} \lim_{R \to + \infty} \int_R^{-R} \frac{p'(i s)}{p(i s)} ds = \frac 1{2 \pi} \lim_{R \to + \infty} \int_R^{-R} \frac{p'_i(s) - i p'_r(s)}{p_r(s) + i p_i(s)} ds =  \\
= \frac 1{2 \pi} \lim_{R \to + \infty} \int_R^{-R} \frac{(p'_i(s)p_r(s) - p'_r(s)p_i(s)) - i (p'_r(s) p_r(s) + p'_i(s) p_i(s))}{p^2_r(s)+ p^2_i(s)} ds = \\
= \frac 1{2 \pi} \lim_{R \to + \infty} \int_R^{-R} \left[\arg(p_r(s) + i p_i(s)) - \frac 12 i\log(p^2_r(s) + p^2_i(s))\right]'ds =  - \wind(p_r + i p_i).
\end{multline*}
where $p(i s) = p_r(s) + i p_i(s)$, the two polynomials $p_r, p_i$ are related to the $p^r,p^i$ introduced in \eqref{pr pi} by the relation $p_r = p^r$, $p_i = - s p^i$, and $ \wind(p_r + i p_i)$ is the winding number around zero of the curve parameterised by
\[
p_r+ i p_i : \mathbb R \to \mathbb C, \qquad s \mapsto p_r(s) + i p_i(s) = p(i s).
\]
A similar argument can be used in the other half-plane and it gives
\[
2 \beta + \delta - \frac n2  = \wind(p_r + i p_i).
\]
Observe that, by simple considerations on the asymptotic of the curve $s \mapsto p_r(s) + i p_i(s)$, the winding number is an integer if $m$ is even and is a semi-integer if $m$ is odd. Moreover, the derivative of $\arg(p_r + i p_i)$ is a rational function with numerator a polynomial in $s^2$ of degree at most $2m-2$ and denominator a polynomial in $s^2$ of degree $2 m$. The integral can be in turns computed using the residue theorem once again. It follows that
\begin{thm}\label{indices}
With the notations of the previous sections $\gamma$ and $\delta$ can be computed using Sturm's theorem, and
\[
\alpha = \frac12 \left( \frac n2 - \wind(p_r + i p_i) - \gamma \right), \qquad \beta  = \frac 12 \left( \frac n2 + \wind(p_r + i p_i) - \delta \right).
\]
\end{thm}

Let us explicitly write the formulas of Theorem~\ref{indices} in the low-degree cases. We only compute $\gamma$ and $\alpha$, since the other two indices $\delta$ and $\beta$ lead to similar expressions. When $n=2$ we have
\[
\gamma = \var\left(d_2,-d_1, d_1^2-4 d_2\right) - \var\left(1,d_1^2-4 d_2\right)
\]
and
\[
\wind = \frac{1}{2\pi} \int_{-\infty}^{+ \infty} \frac{- d_1 \mu ^2- d_1 d_2}{\mu ^4 + (d_1^2 -2 d_2) \mu ^2 + d_2^2} d\mu.
\]

When $n=3$ we have
\begin{multline*}
\gamma = \var\left(d_3,-d_2,d_1 d_2- 9 d_3,4 d_3 d_1^3-d_2^2 d_1^2-18 d_2 d_3 d_1+4 d_2^3+27 d_3^2\right) + \\
- \var\left(-1,3 d_2- d_1^2,4 d_3 d_1^3-d_2^2 d_1^2-18 d_2 d_3 d_1+4 d_2^3+27 d_3^2\right),
\end{multline*}
\[
\wind = \frac 1{2\pi} \int_{-\infty}^{+ \infty} \frac{-d_1 \mu ^4 + (3 d_3 - d_1 d_2) \mu ^2 - d_2 d_3}{\mu ^6 + (d_1^2 -2 d_2) \mu ^4+ (d_2^2 -2 d_1 d_3)\mu ^2 + d_3^2} d\mu.
\]
When $n=4$ the expressions for $\gamma$ become cumbersome, while
\[
\wind = \frac 1{2\pi} \int_{-\infty}^{+ \infty} \frac{-d_1 \mu ^6 + (3 d_3 - d_1 d_2) \mu ^4 + (3 d_1 d_4 - d_2 d_3) \mu ^2 - d_3 d_4}{\mu ^8 + (d_1^2  - 2 d_2) \mu ^6 + (d_2^2 -2 d_1 d_3 + 2 d_4) \mu ^4 + (d_3^2 - 2 d_2 d_4) \mu ^2+d_4^2} d\mu.
\]

\section{Application}\label{application}

To conclude, we want to show how the algebraic equations that define $\mathcal Z$, $\mathcal D$, and $\mathcal R$ become useful in an application. Let us consider the celebrated Lorenz system, that is the vector field
\[
X = \begin{pmatrix} a (y - x) \\ b x - y - x z \\ x y - c z \end{pmatrix}.
\]
The equilibria of $X$ are $e = (0,0,0)$ and $e_\pm = (\pm \sqrt{c(b-1)}, \pm \sqrt{c(b-1)},b-1)$. Let us consider the equilibrium $e$, which always exists. The linearisation, the characteristic polynomial $p$, and the associated polynomials $q^r,q^i$ at such equilibrium are
\[
JX_e = \begin{pmatrix}  -a & b & 0 \\ a & -1 & 0 \\  0 & 0 & -c \end{pmatrix}, \qquad \begin{cases} p = -\lambda ^3 - \lambda ^2 (1+a+c) + \lambda  (a b-a c-a-c)+ a c (b -1) \\
q^r =  (a+c+1)\nu + a c (b-1)\\
q^i = \nu + a b-a c-a-c.
\end{cases}
\]

\begin{figure}
\begin{center}
\includegraphics[width=4cm]{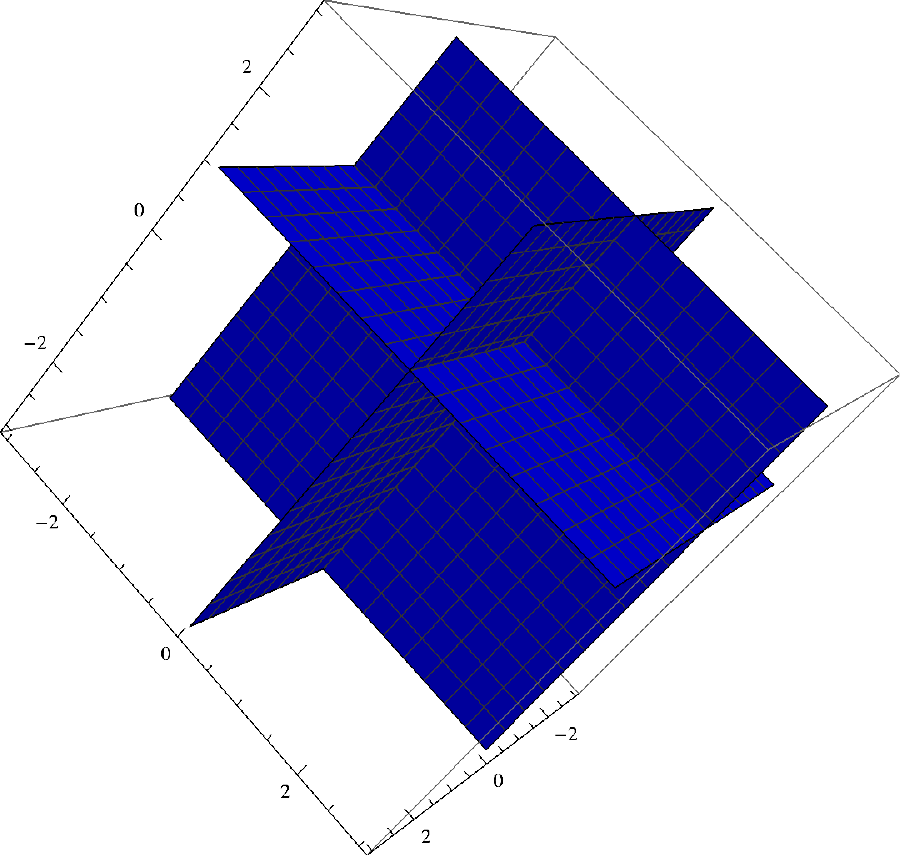}\includegraphics[width=4cm]{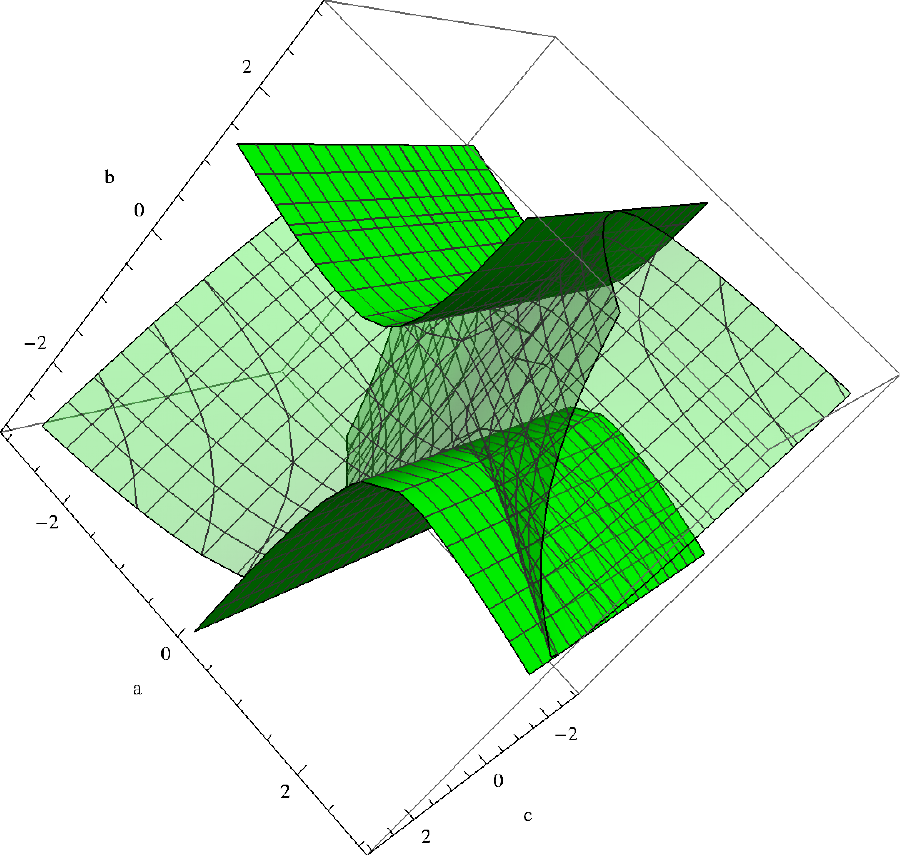}
\includegraphics[width=4cm]{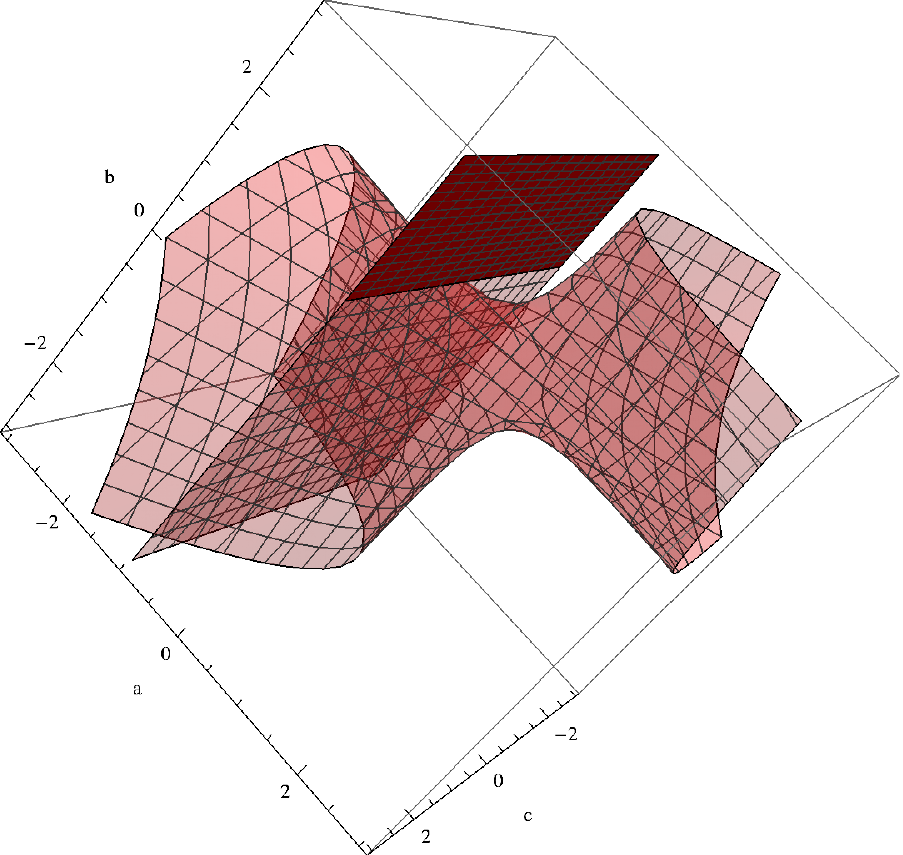}
\end{center}
\caption{In the first pane the determinant locus  $\mathcal Z$, in the second pane the discriminant locus $\mathcal D$ which is the union of two varieties, one across which $\delta$ changes sign (solid green) and another across which $\delta$ does not change sign (transparent green), in the third pane the resultant locus $\widetilde{\mathcal R}$ in solid and transparent red. In solid red its semialgebraic subvariety $\mathcal R$.}\label{Lorentz e0}
\end{figure}

\begin{figure}
\begin{center}
\includegraphics[width=7cm]{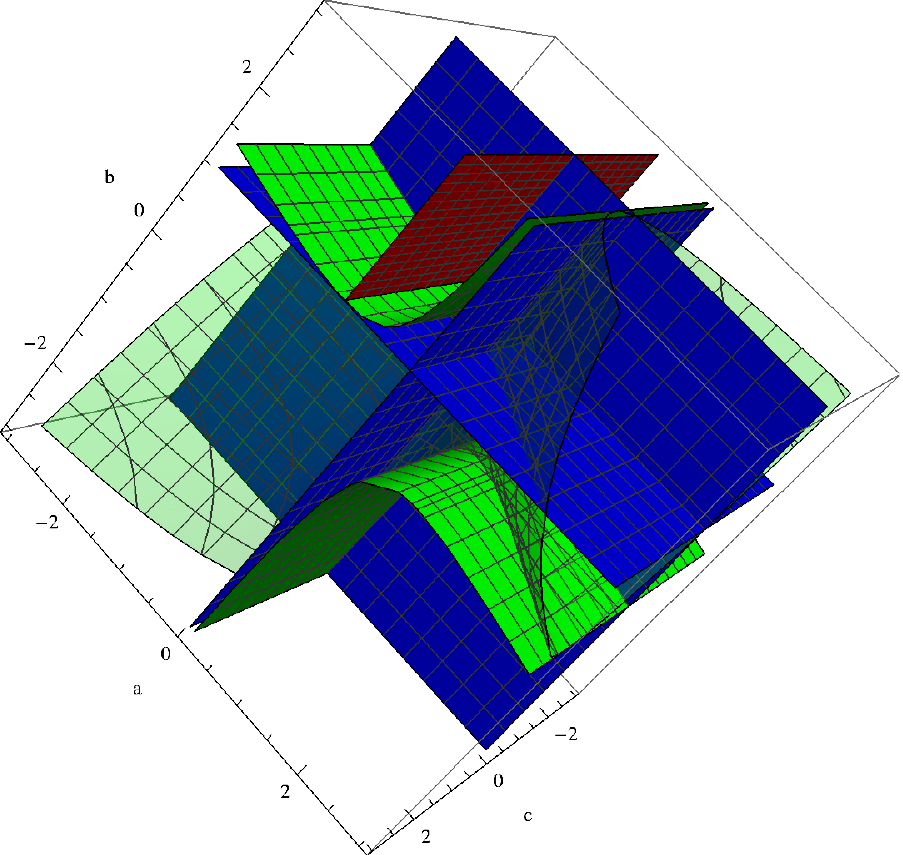} \qquad \includegraphics[width=7cm]{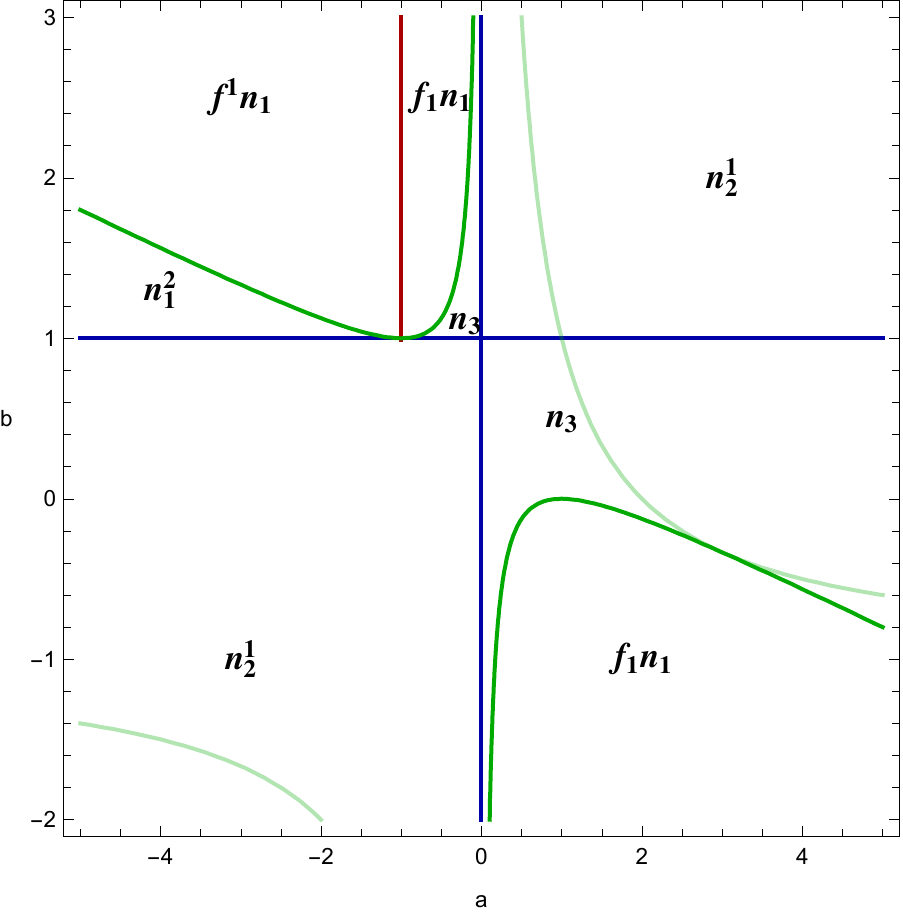}
\end{center}
\caption{In the first pane all the relevant marginal loci together, in the second pane is a section at $c = 2$ of the three-dimensional plot with the indication of the spectral types of the equilibrium $e_0$. The transparent green curve is a discriminant locus in which two negative real eigenvalues touch but on both sides of the curve they separate in the real axis.}\label{Lorentz e0 bis}
\end{figure}

In parameter space $(a,b,c)$ we have that
\[
\zeta = \zeta_1 \, \zeta_2 \, \zeta_3, \text{ with } \quad \zeta_1 = a, \quad \zeta_2 = b-1, \quad \zeta_3 = c,
\]
\[
\rho = \rho_1 \rho_2, \text{ with } \quad \rho_1 = 1+a, \quad \rho_2 = a - a b + c + a c + c^2, \quad \text{and } \quad \sigma = a - a b + c + a c,
\]
\[
\delta = \delta_1 \delta_2^2, \text{ with } \quad \delta_1 = (-1 + a)^2 + 4 a b, \quad \delta_2 = (-1 + c) c - a (-1 + b + c).
\]
This allows to draw the marginal loci in Figure~\ref{Lorentz e0}. The cumulative picture of the marginal locus and a section at $c=2$ with the indication of the spectral types is shown in Figure~\ref{Lorentz e0 bis}.

 \section{Conclusions}

In this article we used a few fundamentals facts:
\begin{itemize}
\item the marginal loci are stratified, their regular stratum corresponds to a simple degeneracy. With simple degeneracy we mean that either zero is a simple root of $p$ ($\mathcal Z$), or $p$ has a double root in the real line ($\mathcal D$), or $p$ has two roots in the imaginary axis ($\mathcal R$). The strata with higher codimension correspond to higher degeneracies (e.g. three roots that coincide in the real axis, two couples of complex conjugate roots that coincide);
\item generically, across the regular stratum of the marginal loci, the corresponding vanishing function, $\zeta$, $\delta$, and $\rho$ change sign, and respectively one root of $p$ changes sign crossing $\mathcal Z$, two roots of $p$ change from two real to a couple of conjugate numbers crossing $\mathcal D$, two complex conjugate roots of $p$ have real part that changes sign crossing $\mathcal R$;
\item if, crossing a marginal locus, the corresponding vanishing function does not change sign, then the transverse generic change of spectral type of the previous point does not take place.
\item the penultimate Euclid's remainder of two polynomials $p,q$ is generically a polynomial of degree one whose zero is the unique common zero of the two polynomials when the ultimate Euclid's remainder, which is the resultant of $p,q$, vanishes (if $q = p'$ the same is true with double root of $p$ instead of unique common zero and discriminant instead of resultant);
\end{itemize}

These facts are basic notions of singularity theory and algebraic geometry, and we think that giving formal justifications in this article would obscure the relevant information given here. The applications of these ideas to relevant parametric-dependent systems will prove extremely useful.

We found enlightening to numerically draw the loci in parameter space, numerically compute the location of zeroes of the characteristic polynomial in the complex plane, numerically compute the indices $\alpha,\beta,\gamma,\delta$ through the variations of Sturm's sequences and the winding number given in Section~\ref{7}, and visualise the changes of spectral type with a Mathematica manipulation.

\section*{Acknowledgments}
The author acknowledges the financial support of Università degli Studi di Catania, progetto PIACERI \emph{Analisi qualitativa e quantitativa per sistemi dinamici finito e infinito dimensionali con applicazioni a biomatematica, meccanica, e termodinamica estesa classica e quantistica}, PRIN 2017YBKNCE \emph{Multiscale phenomena in Continuum Mechanics: singular limits, off-equilibrium and transitions}, and GNFM (INdAM). I also thank Francesco Russo for fruitful conversations.

\printbibliography[heading=bibliography]
%

\end{document}